%% file: demedicis_stanton_white.tex

\input amstex


\input texdraw

\def\btableaux#1{
\btexdraw		
\ifx#1r \def\dx{-1} \else \def\dx{1} \fi	
\def\shading{0.65}	
\setunitscale{.5}	
\drawdim cm		
\linewd 0.02	   	
\textref h:C v:C	
\countdef\bydir=11	
\countdef\bxdir=12	
\bydir=0
\bxdir=0
}

\def\etableaux{
\etexdraw		
}

%
%

\def\plbox#1{\bbox#1{n}}

\def\shbox#1{\bbox#1{y}} 

\def\nbox{
\bxdir=0
\advance\bydir by -1
}

\def\bbox#1#2{
\move ({\the\bxdir} {\the\bydir})
\rlvec (1  0)
\rlvec (0  1)
\rlvec (-1 0)
\rlvec (0 -1)
\ifx#2y \lfill f:{\shading} \fi

\rmove (.5 .5)
\htext {#1}
\rmove (-.5 -.5)
\advance\bxdir by \dx
} 

%
%

\documentstyle{amsppt}
\define\qba{\left[\matrix n_2\\n_2+n_3-n_1\endmatrix\right]_q}
\define\qbb{\left[\matrix n_3\\n_2+n_3-n_1\endmatrix\right]_q}
\font\gothic=eufm10
\define\GS{\text{\gothic S}}
\define\lo{\underline}
\define\hi{\overline}
\magnification=1200
\baselineskip=12pt
\topmatter
\NoRunningHeads
\title
The Combinatorics of $q$-Charlier Polynomials
\endtitle
\author
Anne de M\'edicis\footnote"*"{\hbox{This work was supported by NSERC funds.}\hfill}
\\
Dennis Stanton\footnote"$\dag$"{This work was supported by NSF grant DMS-9001195.\hfill\break}\\
Dennis White
\endauthor
\address{School of Mathematics,
University of Minnesota, Minneapolis, MN 55455.}
\endaddress
\abstract{We describe various aspects of the Al-Salam-Carlitz 
$q$-Charlier polynomials.
These include combinatorial descriptions of the moments, 
the orthogonality relation, and the linearization coefficients.}
\endabstract
\endtopmatter
\NoBlackBoxes
\nologo
\document
\baselineskip=12pt


\subheading{1. Introduction}

The Charlier polynomials $C_{n}^a(x)$ are well-known analytically 
\cite{4}, and have been studied combinatorially by various authors
 \cite{8}, \cite{12}, \cite{16}, \cite{17}, \cite{20}. 
The moments for the measure of these orthogonal polynomials are
$$
\mu_n=\sum_{k=1}^n S(n,k)a^k,
\tag1.1
$$
where $S(n,k)$ are the Stirling numbers of the second kind. The purpose of 
this paper is to study combinatorially an appropriate $q$-analogue of
$C_{n}^a(x)$, whose moments are  a $q$-Stirling version of (1.1). 
While  studying these polynomials, we 
use statistics on set partitions which are $q$-Stirling distributed.

Our main result (Theorem 3) is the combinatorial proof of the linearization
coefficients for these polynomials.
In the $q=1$ case, the linearization coefficients are given as a
polynomial in $a$, whose coefficients are quotients of factorials
(see (4.4)). This has a simple combinatorial explanation. However,
in the $q$-case the coefficients are not the analogous quotients of
$q$-factorials. They are alternating sums of quotients
of $q$-factorials, and thus a combinatorial explanation is much more
difficult. 
From the combinatorial interpretations of the polynomials and their
moments, in terms of weighted partial permutations and set partitions,
we deduce a combinatorial interpretation for the linearization 
coefficients of a product of three $q$-Charlier polynomials. We then
apply a weight-preserving sign-reversing involution defined in five steps.
Theorem 3 is obtained by enumerating the remaining fixed points. Some
of the steps of the involution are quite straight-forward, but some others
are more complicated. They use more sophisticated techniques such as
encoding of permutations or set partitions into {\it 0--1 tableaux}
(cf \cite{6}, \cite{18}), which are fillings of Ferrers diagrams with
0's and 1's such that there is exactly one 1 in each column.
They also use {\it interpolating statistics} on set partitions, as were
introduced by White in \cite{22}. Indeed, the characterization of the 
final set of fixed points uses a bijection $\Psi_S$ of White between
interpolating statistics, making their enumeration all the more complicated.

It turns out that our $q$-Charlier polynomials 
are not what have classically been called
$q$-Charlier; in fact they are rescaled versions of the Al-Salam-Carlitz 
polynomials \cite{4, p.196}. Some comparisons to the classical $q$-Charlier
are given in \S7. Zeng \cite{24} has also studied both families of
 polynomials from the
associated continued fractions.

The basic combinatorial interpretation of the polynomials is 
given in Theorem 1. Several facts about the polynomials can be proven 
combinatorially. The combinatorics of set partitions, restricted growth
functions and 0--1 tableaux is discussed in \S3, and
the statistic for the moments is given in Theorem 2. In \S4, we
state our main theorem, Theorem 3, giving the linearization coefficient
for a product of three $q$-Charlier polynomials, and we set up the general
combinatorial context for its demonstration. The five steps of the
weight-preserving sign-reversing involution proving Theorem 3
are given in \S5, and the combinatorial evaluation of the remaining
fixed points is the subject of \S6.

We use the standard notation for $q$-binomial coefficients and shifted factorials found in \cite{11}. We will also need
$$
[n]_q =\frac{1-q^n}{1-q},
$$
and
$$
[n]!_q=[n]_q[n-1]_q\cdots [1]_q.
$$


\subheading{2. The $q$-Charlier polynomials}

We define the $q$-Charlier polynomials by the three term recurrence relation
$$
C_{n+1}(x,a;q)=(x-aq^n-[n]_q)C_n(x,a;q)-a[n]_q q^{n-1}C_{n-1}(x,a;q),
\tag2.1
$$
where $C_{-1}(x,a;q)=0$ and $C_{0}(x,a;q)=1$.

It is not hard to show that these polynomials are rescaled versions of the Al-Salam Carlitz polynomials \cite{4, p.196}
$$
C_n(x,a;q)=a^n U_n(\frac{x}{a}-\frac{1}{a(1-q)},\frac{-1}{a(1-q)}).
\tag2.2
$$

Since the generating function of the $U_n(x,b)$ is known \cite{4},
we see that
$$
\sum_{n=0}^{\infty}C_n(x,a;q)\frac{t^n}{(q)_n}=
\frac{(at)_{\infty}(-\frac{t}{1-q})_{\infty}}{(t(x-\frac{1}{1-q}))_{\infty}}.
\tag2.3
$$
This gives the explicit formula
$$
C_n(x,a;q)=\sum_{k=0}^n\left[\matrix n\\k\endmatrix\right]_q
(-a)^{n-k}q^{\binom{n-k}2}
\prod_{i=0}^{k-1}(x-[i]_q).
\tag2.4
$$

Clearly, we want a $q$-version of \cite{16}, which gives 
the Charlier polynomials as a generating function of weighted {\it
partial permutations}, i.e. pairs 
$(B,\sigma)$, where 
$B\subseteq \{1,2,\cdots,n\}=[n]$, and $\sigma$ is a permutation on $[n]-B$.
Thus we need only interpret the individual terms in (2.4) for a combinatorial
interpretation. The inside product can be expanded in terms of the 
$q$-Stirling numbers of the first kind. We let $cyc(\sigma)$ be the 
number of cycles 
of a permutation $\sigma$ and $inv(\sigma)$ be the number of inversions of
$\sigma$ written as a product of disjoint cycles (increasing minima, 
minima first in a cycle).
$$
\prod_{i=0}^{k-1}(x-[i]_q)=\sum_{\sigma\in \GS_k} 
(-1)^{k-cyc(\sigma)} q^{inv(\sigma)} x^{cyc(\sigma)}. 
$$
For the sum over $k$ in (2.4), we sum over all $(n-k)$ subsets 
$B\subseteq [n]$. Let 
$$
inv(B)={\sum_{b\in B} (b-1)},
$$
so that the generating function for these subsets is
$$
\left[\matrix n\\k\endmatrix\right]_q q^{\binom{n-k}2}.
$$
We have established the following theorem. 

\proclaim{Theorem 1} The $q$-Charlier polynomials are given by
$$
\align
C_n(x,a;q)&= \sum_{B\subseteq [n]}\sum_{\sigma\in \GS_{n-B}}
q^{inv(\sigma)+inv(B)}(-1)^{n-cyc(\sigma)}a^{|B|}x^{cyc(\sigma)},\\
&= \sum_{B\subseteq [n]}\sum_{\sigma\in \GS_{n-B}} \omega_q(B,\sigma)
x^{cyc(\sigma)}.\\
\endalign
$$
\endproclaim 

A combinatorial proof of the three-term recurrence relation 
(2.1) can be given using Theorem 1.
An involution is necessary. For more details, we refer the reader to \cite{5}.


\subheading{3. The moments}

An explicit measure for the $q$-Charlier polynomials is known, 
\cite{4, p.196}.
It is not hard to find the $n^{th}$ moment of this measure explicitly.
The result is a perfect $q$-analogue of (1.1)
$$
\mu_n=\sum_{k=1}^{n}S_q(n,k)a^k,
\tag3.1
$$
where $S_q(n,k)$ is the $q$-Stirling number of the second kind, 
given by the recurrence
$$
S_q(n,k)=S_q(n-1,k-1)+[k]_qS_q(n-1,k),
\tag3.2
$$
where $S_q(0,k)=\delta_{0,k}$.
In fact, one sees that \cite{13}
$$
S_q(n,k)=\frac{1}{(1-q)^{n-k}}\sum_{j=0}^{n-k}\binom{n}{k+j}
\left[\matrix k+j\\j\endmatrix\right]_q
(-1)^j.
\tag3.3
$$ 

Clearly (3.1) suggests that there is some statistic on set partitions, whose 
generating function is $\mu_n$. This statistic, $rs$, arises from the 
Viennot theory 
of Motzkin paths associated with the three-term recurrence (2.1) \cite{20}. 
We do not give the details of the construction here. 

However, let us  review some combinatorial facts about 
$q$-Stirling numbers.
Set partitions of $[n]=\{1,2,\ldots,n\}$ 
can be encoded as {\it restricted 
growth functions} (or {\it RG-functions}) as follow: if the blocks of $\pi$
are ordered by increasing minima, the RG-function $w=w_1 w_2 \ldots w_n$ is
the word such that $w_i$ is the block where $i$ is located.  For example,
if $\pi= 1 4 7|2 8|3|5 6 9$, $w=123144124$. Note that set partitions on any
set $A$ can be encoded as RG-functions as long as $A$ is a totally ordered
set.

In \cite{21}, Wachs and White investigated four natural statistics
on set partitions, called $ls$, $lb$, $rs$ and $rb$. They are defined as
follow:
$$
\align
ls(\pi)&=ls(w)=\sum_{i=1}^n |\{j: j<w_i, j 
\text{ appears to the left of position } i\}|,\\
lb(\pi)&=lb(w)=\sum_{i=1}^n |\{j: j>w_i, j 
\text{ appears to the left of position } i\}|,\\
rs(\pi)&=rs(w)=\sum_{i=1}^n |\{j: j<w_i, j 
\text{ appears to the right of position } i\}|,\\
rb(\pi)&=ls(w)=\sum_{i=1}^n |\{j: j>w_i, j 
\text{ appears to the right of position } i\}|.\\
\endalign
$$
Thus in the example, $ls(\pi)=13$, $lb(\pi)=7$, $rs(\pi)=7$ and $rb(\pi)=11$.
They showed, using combinatorial methods, that each had the same distribution
(up to a constant) on the set $RG(n,k)$ of all restricted growth functions
of length $n$ and maximum $k$, and that their generating function was indeed
$S_q(n,k)$ for $rs$ and $lb$
 (respectively $q^{\binom{k}{2}}S_q(n,k)$ for $ls$ and $rb$).

We also use  another encoding of set partitions in terms
of {\it 0--1 tableaux}. A 0--1 tableau is a pair $\varphi=(\lambda,f)$ where
$\lambda=(\lambda_1\geq\lambda_2\geq\ldots\geq\lambda_k)$ is a partition of
an integer $m=|\lambda|$ and $f=(f_{ij})_{1\leq j\leq \lambda_i}$ is a 
``filling'' of the corresponding Ferrers diagram of shape $\lambda$ with
0's and 1's such that there is exactly one 1 in each column. 
0--1 tableaux were introduced
by Leroux in \cite{18} to establish a $q$-log concavity result conjectured
by Butler \cite{3} for Stirling numbers of the second kind.

There is a natural correspondence between set partitions $\pi$ of $[n]$ with
$k$ blocks and 0--1 tableaux with $n-k$ columns of length less than or 
equal to $k$. Simply write the RG-function $w=w_1 w_2\ldots w_n$ associated
to $\pi$ as a $k\times n$ matrix, with a 1 in position $(i,j)$ if $w_j=i$,
and 0 elsewhere. The resulting matrix is row-reduced echelon, of rank $k$, 
with exactly one 1 in each column. A 0--1 tableau (in the third quadrant)
is then obtained by removing all the pivot columns and the 0's that lie
on the left of a 1 on a pivot column. Figure 1 illustrates these 
manipulations for $\pi=1 2 4 7|3 9 (12)|5 6 8 (11)|(10)$.

\midinsert \baselineskip=10pt
$$\vbox{\offinterlineskip
\halign{
&\vrule#&\hfil$\displaystyle{#}$\hfil&\vrule#\cr
\noalign{\hrule}\cr
height 6pt&\omit&\cr
&\ \pmatrix 1 & 1 & 0 & 1 & 0 & 0 & 1 & 0 & 0 & 0 & 0 & 0\\
 0 & 0 & 1 & 0 & 0 & 0 & 0 & 0 & 1 & 0 & 0 & 1\\
 0 & 0 & 0 & 0 & 1 & 1 & 0 & 1 & 0 & 0 & 1 & 0\\
 0 & 0 & 0 & 0 & 0 & 0 & 0 & 0 & 0 & 1 & 0 & 0\\
\endpmatrix&\cr
height 6pt&\omit&\cr
& \Updownarrow&\cr
height 6pt&\omit&\cr
&\ \vcenter{
\vbox{
\btableaux{r}
\plbox{0}\plbox{0}\plbox{0}\plbox{0}\plbox{1}\plbox{0}\plbox{1}\plbox{1}\nbox
\plbox{1}\plbox{0}\plbox{1}\plbox{0}\plbox{0}\plbox{0}\plbox{0}\nbox
\plbox{0}\plbox{1}\plbox{0}\plbox{1}\plbox{0}\plbox{1}\nbox
\plbox{0}\plbox{0}\nbox
\etableaux}
}&\cr height 6pt&\omit&\cr
\noalign{\hrule}}}$$
\baselineskip=10pt
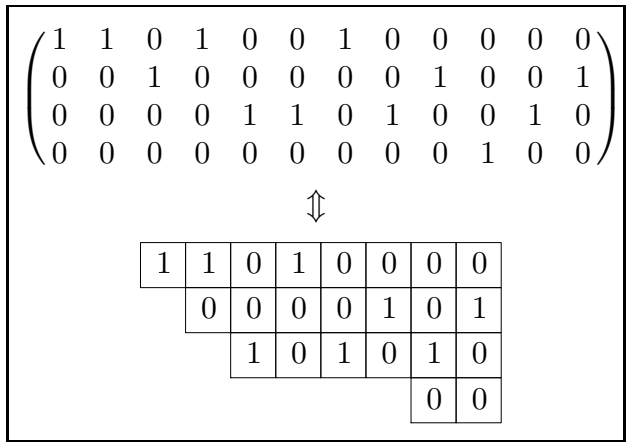
\centerline{Figure 1: Correspondence between partitions and 0--1 tableaux}
\endinsert

We define two statistics on 0--1 tableaux $\varphi$: first, the {\it inversion
number}, $inv(\varphi)$, which is equal to the number of 0's below a 1 in
$\varphi$; and the {\it non-inversion number}, $nin(\varphi)$, which is equal
to the number of 0's above a 1 in $\varphi$. For example, for $\varphi$ in
Figure 1, $inv(\varphi)=7$ and $nin(\varphi)=8$. 
Note that an easy involution
on the columns of 0--1 tableaux sends the inversion number to the non-inversion
number and vice-versa. We  call this map the {\it symmetry involution}.

It is not hard to see that the inversion number (respectively non-inversion
number) on 0--1 tableaux corresponds to the statistic $lb$ (resp. 
$ls-\binom{k}{2}$) on set partitions. 

Similarly, permutations $\sigma$ of $[n]$ in $k$ cycles can  be 
encoded as 0--1 
tableaux with $n-k$ columns of distinct lengths less than or equal to
$n-1$. The correspondence is defined by recurrence on $n$. Suppose $\sigma$
is written as a standard product of cycles. If $n=1$, then $\sigma=(1)$
corresponds to the empty 0--1 tableau $\varphi=\varnothing$. Otherwise, 
let $\sigma\in\GS_{n+1}$ and let
$\varphi$ denote the 0--1 tableau
associated to the permutation $\sigma$ in which
 $(n+1)$ has been erased. There are two cases. If $(n+1)$ is the
minimum of a cycle in $\sigma$, then $\sigma$
corresponds to $\varphi$. If $(n+1)$ is not the
minimum of a cycle, then it appears in $\sigma$ at a certain 
position $i$, $2\leq
i\leq n+1$. The permutation $\sigma$ then corresponds to the 0--1 tableau
$\varphi$ plus a column of length 
$n$ with a 1 in the ($i-1$)-th position (from top to bottom). For example,
$\sigma=(1,3,4,7,2)(5,6)(8)$ corresponds to the following 0--1 tableau.
\midinsert \baselineskip=10pt
\centerline{\vbox{
\btableaux{l}
\plbox{0}\plbox{0}\plbox{0}\plbox{1}\plbox{1}\nbox
\plbox{0}\plbox{0}\plbox{1}\plbox{0}\nbox
\plbox{1}\plbox{0}\plbox{0}\nbox
\plbox{0}\plbox{0}\nbox
\plbox{0}\plbox{1}\nbox
\plbox{0}\nbox
\etableaux}}
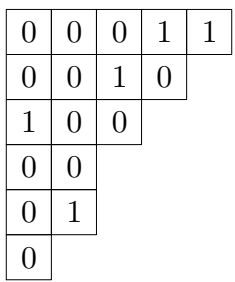
\centerline{Figure 2: Correspondence between permutations and 0--1 tableaux}
\endinsert
It is not hard to see that under this transformation,
the inversion number on 0--1 tableaux
corresponds to the inversion number on permutations, as defined in
\S2. Thus, their generating functions are the $q$-Stirling numbers of 
the first kind $c_q(n,k)$.

In \cite{6}, de M\'edicis and Leroux investigated $q$ and $p,q$-Stirling
numbers from the point of view of the unified 0--1 tableau approach. In
particular, they proved combinatorially or algebraically 
a number of identities
involving $q$-Stirling numbers.

For the combinatorial interpretation of the moments of the
$q$-Charlier polynomials in terms of set partitions 
$\pi$, we need two statistics. The number of blocks $\#blocks(\pi)$ is one,
and the other statistic is $rs(\pi)$.

\proclaim{Theorem 2} The $n^{th}$ moment for the $q$-Charlier polynomials
is given by
$$
\mu_n=\sum_{\pi\in P(n)} a^{\#blocks(\pi)} q^{rs(\pi)}.
$$
\endproclaim

As we mentioned,
many other $q$-Stirling distributed statistics have been found \cite{21}. 
It is surprising that the Viennot theory naturally gives a so-called ``hard"
statistic ($rs$), not an easy one (e.g. $lb$, \cite{21}). 
Other variations on the 
$rs$-statistic 
can be given from the Motzkin paths, although the 
$lb$-statistic 
 is not among them. It can be derived from the Motzkin paths 
associated with the ``odd" polynomials for (2.1). 


\subheading{4. The orthogonality relation and the linearization of products}

Let $L$ be the linear functional on polynomials that corresponds to
integrating with respect to the measure for the Charlier polynomials. 
The orthogonality relation is
$$
L(C_{n}^a(x)C_{m}^a(x))= a^nn!\ \delta_{m,n}.
\tag4.1
$$
The $q$-version of (4.1) is
$$
L_q(C_{n}(x,a;q)C_{m}(x,a;q))=a^nq^{\binom{n}{2}}[n]!_q \delta_{m,n}.
\tag4.2
$$
Since the polynomials $C_{n}(x,a;q)$ and $L_q$ have combinatorial
definitions from Theorems 1 and 2, it is possible to restate (4.2) as a
combinatorial problem. We will give an involution which then proves (4.2)
in this framework.

A more general question is to find
$L(C_{n_1}^a(x)C_{n_2}^a(x)\cdots C_{n_k}^a(x))$
for any $k$. A solution is equivalent to finding the coefficients
$a_{n_k}$ in the expansion
$$
C_{n_1}^a(x)C_{n_2}^a(x)\cdots C_{n_{k-1}}^a(x) = \sum_{n_k} a_{n_k} C_{n_k}^a(\
x).
$$
This had been done bijectively for some classes of Sheffer orthogonal
polynomials in \cite{5}, \cite{7}, \cite{9}, \cite{10}. Moreover,
in the $q$-case of Hermite polynomials, some remarkable consequences have been
found \cite{15}.

For the Charlier polynomials, it is easy to see that
$$
\sum_{n_1,\cdots,n_k=0}^{\infty}
L(C_{n_1}^a(x)C_{n_2}^a(x)\cdots C_{n_k}^a(x))
\frac{t_1^{n_1}}{n_1!}\cdots
\frac{t_k^{n_k}}{n_k!}= e^{a(e_2(t_1,\cdots,t_k)+\cdots+e_k(t_1,\cdots,t_k))},
\tag4.3
$$
where $e_i$ is the elementary symmetric function of degree $i$, \cite{19}.
In this case  $L(C_{n_1}^aC_{n_2}^a\cdots C_{n_k}^a)$
is a polynomial in $a$ with
positive integer coefficients;  a combinatorial interpretation
of this coefficient has been given (\cite{12} and \cite{23}).
For $k=3$, (4.3) is equivalent to
$$
L(C_{n_1}^a(x)C_{n_2}^a(x) C_{n_3}^a(x))=
\sum_{l=0}^{\lfloor(n_1+n_2-n_3)/2\rfloor}{a^{n_3+l}n_1!n_2!n_3!
\over l! (n_3-n_2+l)!  (n_3-n_1+l)! (n_1+n_2-n_3-2l)!}.
\tag 4.4
$$

One can  hope that
$L_q(C_{n_1}(x)C_{n_2}(x) C_{n_3}(x))$
is simply a weighted version, with an appropriate
statistic, of the $q=1$ case. However this is false. For example,
$$
L_q(C_{2}(x)C_{2}(x)C_{1}(x))=q(q^2+2q+1)a^2+q(q^3+q^2-q-1)a^3.
$$
Nonetheless,
we have an exact formula for $L_q(C_{n_1}(x,a,q)C_{n_2}(x,a,q)C_{n_3}(x,a,q))$,
which is equivalent to one of Al-Salam-Verma \cite{1}.

\proclaim{Theorem 3} Let $n_3\geq n_1\geq n_2\geq 0$. Then
$$
\align
L_q(C_{n_1}(x)&C_{n_2}(x)C_{n_3}(x))=
\sum_{l=0}^{n_1+n_2-n_3}\sum_{j=0}^la^{n_3+l}
q^K(q-1)^{l-j}{[n_1-j]!_q\over
[n_1-l]!_q}\left[\matrix n_2\\l-j\endmatrix\right]_q\\ &
[n_3]!_q\left[\matrix n_1\\j\endmatrix\right]_q
\left[\matrix n_2-l+j\\n_3-n_1+j\endmatrix\right]_q
{[j]!_q[n_1-j]!_q\over [n_3-n_2+l]!_q}
\left[\matrix n_1+n_2-n_3-l\\j\endmatrix\right]_q,\tag 4.5\\
\endalign
$$
where
$$
\align
K=&\binom{l-j}{2}+\binom{n_1}{2}+j(-n_3-j+1)+\binom{j}{2}
+\binom{n_2-l+j}{2}\\&+(n_3-n_1+j)(n_3-n_2+l)
+j(n_3-n_2+l).\\
\endalign
$$
\endproclaim

The generating function of $L_q(C_{n_1}(x)C_{n_2}(x)C_{n_3}(x))$ can
be evaluated from Theorem 3, yielding
$$
\align
\sum_{n_1,n_2,n_3}L_q(C_{n_1}(x)&C_{n_2}(x)C_{n_3}(x))
{t_1^{n_1}\over [n_1]!_q}{t_2^{n_2}\over [n_2]!_q}{t_3^{n_3}\over [n_3]!_q}=\\
&
(-t_3;q)_{\infty}(-at_1t_2(1-q);q)_{\infty}\,_2\phi_1
\left({at_1(1-q),\, at_2(1-q)\atop -at_1t_2(1-q)}; q, -t_3\right).
\tag4.6\\
\endalign
$$
Letting $q\rightarrow 1$ in (4.6) gives back (4.3)
for $k=3$.
This generating function can also be evaluated directly using the measure
(\cite{4, p.196}), the generating function (2.3) for the polynomials  and a
$\,_3\phi_2$ transformation.

More generally, for $k\geq 4$, the generating function of $L_q(C_{n_1}(x)
\ldots C_{n_k}(x))$ can be expressed as a difference of two basic 
hypergeometric series. This has been done by Ismail and Stanton \cite{14}
for the Al-Salam Carlitz polynomials, so an equivalent formula can be deduced
for the $q$-Charlier polynomials using (2.2).

Let us set up the combinatorial context in which Theorem 3 will be proven.
We first introduce notations and conventions that will be used throughout
the proof. Define
$$
\align
L_q(n_1,n_2,n_3)=\{& ((B_i,\sigma_i);\pi)=((B_1,\sigma_1),(B_2,\sigma_2),
(B_3,\sigma_3);\pi) | \\ & (B_i,\sigma_i) \text{ is a partial permutation
on the set } \{ i\}\times [n_i],\\ &
 \text{ and $\pi$ is a partition on the 
cycles of $\sigma_1$, $\sigma_2$ and $\sigma_3$ }\}.\\
\endalign
$$

We will say that an element of the set $\{ i\}\times [n_i]$ is
{\it of color} $i$. When giving examples of elements of 
$L_q(n_1,n_2,n_3)$, to simplify notation, pairs $(1,i)$, $(2,i)$ and
$(3,i)$ will always be denoted $\underline i$, $i$ and $\overline i$
respectively. Thus a typical element of $L_q(8,7,10)$ would be
described in the following way: $B_1=\{\underline 2,{\lo 3}\}$,
$B_2=\varnothing$, $B_3=\{{\hi 5},{\hi 9},{\hi {10}}\}$, and
$\pi=({\lo 1},{\lo 5},{\lo 7})({\lo 8})({\hi 1})|({\lo 4}) |
({\lo 6})(3,5)({\hi 3},{\hi 7})|(1,4,2)(6,7)|(\hi 2,\hi 8,\hi 4,\hi 6)$
(the underlying permutations $\sigma_1=({\lo 1},{\lo 5},{\lo 7})
({\lo 4})({\lo 6})({\lo 8})$, $\sigma_2=(1,4,2)(3,5)(6,7)$ and
$\sigma_3=({\hi 1})({\hi 2},{\hi 8},{\hi 4},{\hi 6})({\hi 3},{\hi 7})$
can be recovered from $\pi$).

Note that the lexicographic order on pairs $(i,j)$ induces a total
order on the cycles of $\sigma_1$, $\sigma_2$ and $\sigma_3$, according to 
their minima. Therefore we can talk about RG-functions. We will always
use the letter $w$ to denote the RG-function associated to $\pi$. 
In the above example, $w=1231434153$. The first $cyc(\sigma_1)$ letters
of $w$ correspond to the positions of cycles of color 1 in $\pi$,
the next  $cyc(\sigma_2)$  to the positions of cycles of color 2,
and the last  $cyc(\sigma_3)$ letters to the positions of
cycles of color 3. We will denote by $w_a$, $w_b$ and $w_c$ respectively
these portions of $w$. In the above example, we have 
$w_a=1231$, $w_b=434$, $w_c=153$, and $w=w_aw_bw_c$, the concatenation
of words $w_a$, $w_b$ and $w_c$.

Finally, we will use the notation $Supp(w)$ (or $Supp(\sigma)$ or
$Supp(\pi_i)$) to denote the underlying set of letters of a word $w$
(or a permutation $\sigma$ or a block $\pi_i$ of a partition $\pi$
respectively).

From Theorems 1 and 2,
we deduce that
$$
L_q(C_{n_1}(x)C_{n_2}(x)C_{n_3}(x))=
\sum_{((B_i,\sigma_i);\pi)\in L_q(n_1,n_2,n_3)}\omega_q((B_i,\sigma_i);\pi),
\tag 4.7
$$
where
$$
\omega_q((B_i,\sigma_i);\pi)=\omega_q(B_1,\sigma_1)\omega_q(B_2,\sigma_2)
\omega_q(B_3,\sigma_3)q^{rs(\pi)}a^{\#blocks(\pi)},\tag 4.8
$$
and $\omega_q(B,\sigma)$ was defined in Theorem 1, as a signed monomial in
the variables $a$ and $q$. This gives a combinatorial
interpretation of the right-hand side of (4.5).

For $q=1$, the negative coefficients of $a$ are
counterbalanced by the positive coefficients of $a$, and (4.7) is a
polynomial with positive coefficients.
Indeed, in that case, it is not hard to find a weight-preserving
sign-reversing involution on $L_q(n_1,n_2,n_3)$ (cf \cite{5})
whose fixed points $((B_i,\sigma_i);\pi)$ are characterized by

\noindent{i)}\enspace\  $B_i=\varnothing$ and $\sigma_i=$ Identity, for
$i=1,2,3$;

\noindent{ii)}\enspace the word
 $w_a$ (respectively $w_b$ and $w_c$) contains all distinct
letters, and $Supp(w_a)\subseteq$\break 
\phantom{ii)\enspace}$Supp(w_bw_c)$ (respectively
$Supp(w_b)\subseteq Supp(w_aw_c)$ and
$Supp(w_c)\subseteq Supp(w_aw_b)$).

Identity (4.4) easily follows from $\omega_1$-counting these fixed points.

However, the general $q$-case is much harder, 
and  some negative weights remain.
The sign of $\omega_q((B_i,\sigma_i);\pi)$ comes from the cardinalities
of the sets $B_i$ and the signs of the permutations $\sigma_i$.
In our proof, we  successively apply five weight-preserving
sign-reversing involutions $\Phi_i$ to $ L_q(n_1,n_2,n_3)$, each one
acting on the fixed points of the preceding one.
$\Phi_1$ forces $\sigma_3=Id$, $\Phi_2$ forces $B_3=\varnothing$,
$\Phi_3$ forces $\sigma_1=Id$, $\Phi_4$ forces $B_1=\varnothing$,
and $\Phi_5$ forces $\sigma_2=Id$, leaving $B_2$ arbitrary. Hence the
negative part of (4.7) is due only to $B_2$.

The final set of fixed points, $Fix\Phi_5$, does not contain the
fixed point set (above) for $q=1$. Instead there is a bijection from a subset
of $Fix\Phi_5$ to this set, but it does not preserve the powers
of $q$.

The five weight-preserving sign-reversing involutions $\Phi_i$ and their
respective fixed points sets $Fix\Phi_i$ are given 
in the next section and
the complete characterization of $Fix\Phi_5$ is given by the conditions 
Fix.1 through Fix.4, stated at the beginning of \S6. 
In \S6, we show that the $\omega_q$-weight of $Fix\Phi_5$ is 
equal to the right-hand side of (4.5), thus establishing Theorem 3.


\subheading{5. The weight-preserving sign-reversing involutions $\Phi_i$}

Let us recall that a {\it weight-preserving sign-reversing involution}
(or {\it WPSR-invo\-lu\-tion}) $\Phi$ with weight function $\omega$
is an involution such that for any $e\not\in Fix\Phi$,
$\omega(\Phi(e))=-\omega(e)$.

\subheading{Involution $\Phi_1$} This WPSR-involution will kill any
$((B_i,\sigma_i);\pi)$ such that $\sigma_3$ is not the identity.

Remember that the cycles of $\sigma_3$ are ordered by increasing minima.
Find the greatest cycle $c_{i_0}$ such that either this cycle is of
length $\geq 2$ or it lies in the same block $\pi_i$ of $\pi$ as some
other 1-cycle greater than it. If $c_{i_0}$ satisfies the latter 
condition, the 1-cycle greater than $c_{i_0}$ in the leftmost block
$\pi_h$ of partition $\pi$ is glued to the end of $c_{i_0}$. Then, if
$h=h_0<h_1<\ldots <h_m=i$ denote the indices of the blocks between $\pi_h$
and $\pi_i$ containing 1-cycles greater than $c_{i_0}$, these 1-cycles
are moved from block $\pi_{h_l}$ to block $\pi_{h_{l-1}}$.

For example, for $((B_i,\sigma_i);\pi)\in L_q(9,0,10)$ such that
$B_1=B_2=\varnothing$, $B_3=\{\hi{10}\}$ and
$\pi=(\lo 1)(\hi 1,\hi 2)(\hi 6)|(\lo 2,\lo 8)(\hi 5)|(\lo 3)(\lo 4)|
(\lo 5)(\hi 3,\hi 9)(\hi 8)|(\lo 6,\lo 9)(\lo 7)(\hi 7)|(\hi 4)$, we have
$\sigma_3=(\hi 1,\hi 2)(\hi 3, \hi 9)(\hi 4)
(\hi 5)(\hi 6)(\hi 7)(\hi 8)$,
$c_{i_0}=(\hi 3, \hi 9)$, and
$\Phi_1((B_i,\sigma_i);\pi)$ is given by the same $B_i$'s, $\sigma_3$
becomes $(\hi 1,\hi 2)(\hi 3, \hi 9, \hi 6)(\hi 4)
(\hi 5)(\hi 7)(\hi 8)$, and
$\pi=(\lo 1)(\hi 1,\hi 2)(\hi 5)|(\lo 2,\lo 8)(\hi 8)|(\lo 3)(\lo 4)|
(\lo 5)(\hi 3,\hi 9,\hi 6)|(\lo 6,\lo 9)(\lo 7)(\hi 7)|(\hi 4)$.

Note that the number of inversions gained in $\sigma_3$ is counterbalanced
by the loss in the statistic $rs(\pi)$. Conversely, if $c_{i_0}$ is of
length $\geq 2$ and does not lie in the same block as any other greater
cycles, its image is defined in the obvious way so that $\Phi_1$ is 
an involution. For more details, see \cite{5}.

\subheading{Fixed points for $\Phi_1$} 
The cycle $c_{i_0}$ is not defined if and only if $\sigma_3$
contains only 1-cycles which all lie in different blocks of $\pi$.
Therefore,
$$
\align
Fix\Phi_1=\{&((B_i,\sigma_i);\pi)\in L_q(n_1,n_2,n_3) \mid
\sigma_3 \text{ is the identity}\\ &\text{ and } w_c \text{ contains all
distinct letters}\}.\\
\endalign
$$

\subheading{Involution $\Phi_2$} This WPSR-involution is designed to 
discard all $((B_i,\sigma_i);\pi)\in Fix\Phi_1$ such that $B_3$ is not
empty. 

Let $((B_i,\sigma_i);\pi)\in Fix\Phi_1$ and let $k=\#blocks(\pi)$.

Denote by $j_0$, $0\leq j_0\leq (n_3-1)$, the integer such that
$\overline {j_0+1}=\hbox{min}(B_3)$. 
If $B_3=\varnothing$, we let $j_0=\infty$.
Likewise, denote by $j_1$, $1\leq j_1\leq n_3$, the maximum integer
such that the 1-cycle $(\overline j_1)$ forms a singleton block in $\pi$.
Remember that $\sigma_3=Id$ and $w_c$ contains all distinct letters.
By maximality, $(\overline j_1)$ lies in the $k$-th block of $\pi$. 
Denote by $j_1'$ its contribution to the statistic $rs$, that is the
number of (different) letters after the only occurrence of $k$ in
$w_c$ (and in $w$). If there are no such singleton blocks in $\pi$,
let $j_1=j_1'=\infty$.

There are two cases: $j_0\leq j_1'$, or $j_0>j_1'$.
If $j_0\leq j_1'$, $\Phi_2((B_i,\sigma_i);\pi)$ is obtained by inserting
the 1-cycle $(\overline {j_0+1})$ 
in $\sigma_3$ and by inserting the letter $(k+1)$ in 
$w_c$ at the $(j_0+1)$-th position from the end of $w_c$, leaving everything
else fixed.

For example, for $((B_i,\sigma_i);\pi)$ defined by
$B_1=\varnothing=B_2$, $B_3=\{ \hi 2, \hi 6,\hi 8\}$ and
$\pi=(\lo 1,\lo 6)(5)(\hi 7)|(\lo 2)(1,3,2)(\hi 4)|(\lo 3,\lo 5,\lo 4)|
(4)(\hi 9)|(\hi 1)|(\hi 3)|(\hi 5)$,
$w_c=562714$, $j_0=1$, $j_1=5$ and $j_1'=2$. 
Then the new $w_c$ in $\Phi_2((B_i,\sigma_i);\pi)$ is $w_c=5627184$, and
$\Phi_2((B_i,\sigma_i);\pi)$ is defined by 
$B_1=\varnothing=B_2$, $B_3=\{  \hi 6,\hi 8\}$ and
$\pi=(\lo 1,\lo 6)(5)(\hi 5)|(\lo 2)(1,3,2)(\hi 3)|(\lo 3,\lo 5,\lo 4)|
(4)(\hi 9)|(\hi 1)|(\hi 2)|(\hi 4)|(\hi 7)$.

Note that $\Phi_2((B_i,\sigma_i);\pi)$ has its $j_1'$ equal to the $j_0$
associated to $((B_i,\sigma_i);\pi)$. Conversely, if $j_1'<j_0$,
the image of  $((B_i,\sigma_i);\pi)$ is defined in the obvious way so that
$\Phi_2$ is an involution. $\Phi_2$ is also weight-preserving and 
sign-reversing. For more details, see \cite{5}.

\subheading{Fixed points for $\Phi_2$}
Fixed points correspond to the case $j_0=j_1'=\infty$. This means
that $B_3=\varnothing$ and there are no singleton blocks in $\pi$
of color 3. Therefore,
$$
Fix\Phi_2=\{ ((B_i,\sigma_i);\pi)\in Fix\Phi_1\mid B_3=\varnothing
\text{ and } Supp(w_c)\subseteq Supp(w_aw_b)\}.
$$

Note that $Supp(w_c)\subseteq Supp(w_aw_b)$ is equivalent to the
condition that the $w_aw_b$ is an RG-function whose maximum equals
$\#blocks(\pi)$.

To do $\Phi_3$ and later $\Phi_5$, we need to describe the contribution
to the statistic $rs$ of the elements of color 1 and 2 in partition $\pi$.
Let $w$ be a word on the alphabet $[k]$. Let $w_{ij}$ denote the subword
of $w$ obtained by discarding letters not equal to $i$ or $j$, $1\leq i<j
\leq k$. For instance, if $w=123144124$, $w_{12}=12112$. Then we can write
$$
rs(w)=\sum_{1\leq i<j\leq k} rs(w_{ij}).
$$

\proclaim{Claim} Let $w$ be an RG-function of maximum $k$
and suppose $w=vv'$. Then
$v$ is an RG-function and
$$
\align
rs(w)&=\sum\Sb 1\leq i<j\leq k, \\{i\notin Supp(v')}\endSb rs(v_{ij})
     +\sum\Sb 1\leq i<j\leq k,\\{i\in Supp(v')}\endSb  ls(v_{ij}) + rs(v')\\
&=:rs(w)|_v + rs(v').\\
\endalign
$$
\endproclaim

Thus the contribution to the statistic $rs(w)$ of the initial word $v$,
$rs(w)|_v$, is indeed an interpolation between the hard statistic $rs$
and the easy statistic $ls$, as was studied by  White in \cite{22}. He 
showed in particular that these specific interpolating statistics were
$q$-Stirling distributed, meaning that their generating functions over 
$RG(n,k)$ are the $q$-Stirling numbers of the second kind $S_q(n,k)$, up
to a power of $q$. He provides a bijection on $RG(n,k)$ such that the 
mixed statistic is sent to the easy statistic $ls$ (up to a constant).
More precisely,

\proclaim{Lemma 4} Let $S=\{ s_1<s_2<\ldots<s_m\}\subseteq [k]$. There is
a bijection $\Psi_S:RG(n,k)\to RG(n,k)$ such that for any $w\in RG(n,k)$,
$$
\sum\Sb{1\leq i<j\leq k,}\\ i\in S\endSb
rs(w_{ij})+\sum\Sb 1\leq i< j\leq k,\\ i\in[k]-S\endSb ls(w_{ij})
=ls(\Psi_S(w))-\sum_{j=1}^m(k-s_j).\tag 5.1
$$
\endproclaim

\demo{Proof} 
Define $\Psi_i:RG(n,k)\to RG(n,k)$, $1\leq i\leq k-1$ as follow:
\item{i)} if $w\in RG(n,k)$ has a letter $i$ to the right of the first
occurrence of $(i+1)$, then the rightmost letter $i$ is switched to $(i+1)$
and any $(i+1)$ to its right is changed to $i$. For example, 
$\Psi_1(111212332122)=111212332211$.
\item{ii)} if $w$ does not have a letter $i$ to the right of the first 
occurrence of $(i+1)$, then all $(i+1)$'s to its right are switched to $i$'s.
For example, $\Psi_1(1112232)=1112131$.

For convenience, we will set 
$\Psi_k:RG(n,k)\to RG(n,k)$ to be the identity. Now, given
$S=\{ s_1<s_2<\ldots<s_m\}\subseteq [k]$, $\Psi_S$ is defined as follow:
$$
\Psi_S=(\Psi_k\circ\Psi_{k-1}\circ\ldots\circ\Psi_{s_1})
     \circ(\Psi_k\circ\Psi_{k-1}\circ\ldots\circ\Psi_{s_2})\circ\ldots
     \circ(\Psi_k\circ\ldots\circ\Psi_{s_m}).
$$
Note that $\Psi_S$ preserves the positions of the first occurrences.
For more details, the reader is referred to \cite{22}.\hfill$\square$
\enddemo

\subheading{Involution $\Phi_3$}
This next involution is designed to kill any element $((B_i,\sigma_i);\pi)$
such that $\sigma_1$ is not the identity. Note that since the interpolating 
statistics on $w_a$ are $q$-Stirling distributed, it reduces to proving the
orthogonality relation
$$
\sum_{k=m}^n(-1)^{n-k}c_q(n,k)S_q(k,m)=\delta_{n,m}.
$$
But this formula was deduced in Proposition 3.1 of \cite{6} from a
weight-preserving sign-reversing involution on appropriate pairs of 0--1
tableaux. The general idea is to  map $\sigma_1$ and 
$w_a$ bijectively into a pair of 0--1 tableaux, 
using $\Psi_S$ defined in the previous
lemma and the correspondences described in \S1. 
Then we can apply
the WPSR-involution, essentially  shifting the
rightmost shortest column from one 0--1 tableau to the other. $\Phi_3(
(B_i,\sigma_i);\pi)$ is then obtained by replacing $\sigma_1$ and $w_a$
by the new decoded pair of 0--1 tableaux.
Involution $\Phi_5$ will use similar ideas.

We need only specify the bijective coding of $(\sigma_1,w_a)$ into a pair
of 0--1 tableaux. Let $((B_i,\sigma_i);\pi)\in Fix\Phi_2$ and let 
$n=n_1-|B_1|$, $k=cyc(\sigma_1)$ and  $m=max(Supp(w_a))$.
\item{i)} For $\sigma_1$, simply use the correspondence described in \S3
to get a 0--1 tableau $\varphi_1$ with $(n-k)$ columns of distinct length 
$\leq (n-1)$. Note that $inv(\sigma_1)=inv(\varphi_1)$.
\item{ii)} For $w_a$, we first want to reduce the interpolating statistic
$rs(w)|_{w_a}$ to the easy statistic $ls(w_a)$. This is done by applying
$\Psi_S$ defined in the previous lemma to $w_a$, for $S=[m]
\setminus Supp(w_bw_c)$. We then use the correspondence described in 
\S3 to get a 0--1 tableau $\varphi_2$ with $(k-m)$ columns of length
$\leq m$. There is one last technicality: the statistic $ls$ is sent to the 
non-inversion statistic on 0--1 tableaux (up to the constant $\binom{m}{2}$),
therefore we will apply to $\varphi_2$ the symmetry involution exchanging
non-inversions and inversions, so that for its image $\tilde\varphi_2$, we have
$$
rs(w)|_{w_a}= inv(\tilde\varphi_2)+\binom{m}{2}-\sum_{i\in S}(m-i).
$$

Note that $m$ is not modified by the WPSR-involution 
applied to  pairs of 0--1 tableaux, thus insuring 
that the overall involution $\Phi_3$ is well-defined (the new $w$ is still an
RG-function) and weight-preserving. It is also sign-reversing. Details
are left to the reader.

\subheading{Fixed points for $\Phi_3$}
At the 0--1 tableau level, the only fixed pair of 0--1 tableaux is 
$(\varnothing,\varnothing)$, because in that case, it is impossible to move 
columns. But this can happen if and only if $(n-k)=(k-m)=0$, and
therefore $n=k=m=n_1-|B_1|$, $\sigma_1$ is the identity
on $[n_1]-B_1$, and $w_a=1 2\ldots (n_1-|B_1|)$. Therefore
$$
Fix\Phi_3=\{ ((B_i,\sigma_i);\pi)\in Fix\Phi_2 \mid \sigma_1 \text{ is the 
identity and } w_a=1 2 \ldots (n_1-|B_1|)\}.
$$

\subheading{Involution $\Phi_4$}
This involution is  the simplest. Its task is to eliminate elements
$((B_i,\sigma_i);\pi)$ such that $B_1\not=\varnothing$.

Let $((B_i,\sigma_i);\pi)\in Fix\Phi_3$ and let $i_0$ be the smallest
integer, $1\leq i_0\leq n_1$, such that either $\underline i_0\in B_1$, or the 
1-cycle $(\underline  i_0)$ forms a singleton block in $\pi$. 
Then if $\underline i_0\in
B_1$, insert it as a 1-cycle in $\sigma_1$ and as a singleton block in
$\pi$, and vice-versa.

For example, if $B_1=\{\underline 2\}, B_2=B_3=\varnothing$, and 
$\pi=(\underline 1) (\overline 1,\overline 3)|(\underline 3)| (1,2)(\overline
2)$, 
then $i_0=2$ and the image of $((B_i,\sigma_i);\pi)$ under $\Phi_4$ is
$B_1=\varnothing,B_2=B_3=\varnothing$, and 
$\pi=(\underline 1) (\overline 1,\overline 3)|(\underline 2)|
(\underline 3)| (1,2)(\overline2)$.  
 Details are left to the reader.

\subheading{Fixed points for $\Phi_4$} 
$$
Fix\Phi_4=\{((B_i,\sigma_i);\pi)\in Fix\Phi_3\mid B_1=\varnothing \text{ and }
  Supp(w_a)= [n_1]\subseteq Supp(w_bw_c)\}.
$$

\subheading{Involution $\Phi_5$}
This final WPSR-involution will annihilate the remaining $((B_i,\sigma_i);\pi)$
such that $\sigma_2$ is not the identity. It is the only one 
using the hypothesis $n_3\geq n_1\geq n_2$. The principle of the involution
is similar to $\Phi_3$: we will reduce the problem to finding an 
involution for the easy statistic $ls$. 

Let $((B_i,\sigma_i);\pi)\in Fix\Phi_4$, and let $\#blocks(\pi)=n_3+s$.
First, encode
$\sigma_2$ as a 0--1 tableau $\varphi$ with $(n_2-|B_2|-cyc(\sigma_2))$ 
columns of
distinct lengths $\leq (n_2-|B_2|-1)$, using the correspondence described in
\S3. Note that $inv(\sigma_2)=inv(\varphi)$ and that the shortest
column of $\varphi$ is of length at most $cyc(\sigma_2)$. 

For $w_b$, we reduce the interpolating statistic $rs(w)|_{w_aw_b}$ to the
easy statistic $ls$ by applying $\Psi_S$ defined in Lemma 4 to $w_aw_b$,
with $S=[n_3+s]-Supp(w_c)$. Note that since $w_a=12\ldots n_1$ and
$\Psi_S$ preserves first occurrences, $\Psi_S(w_aw_b)=w_a\tilde w_b$ for
some word $\tilde w_b=\tilde b_1\tilde b_2\ldots \tilde b_k$.
Note also that we must have $\{ n_1+1,\ldots, n_3+s\}\subseteq 
Supp(\tilde w_b)$ (because $w_aw_b$ has maximum $(n_3+s)$). 

For example, if  $((B_i,\sigma_i);\pi)\in Fix\Phi_4$ is defined by
$B_1=B_2=B_3=\varnothing$, and
$\pi=(\lo 1)(2)| (\lo 2) (3) (\hi 5)|(\lo 3)(\hi 2)|(\lo 4)(\hi 4)|
(\lo 5)(5)(\hi 1)|(1)(4)(\hi 3)$,  we have
$w_a=12345$, $w_b=61265$, $w_c=53642$, and $\sigma_2=(1)(2)(3)(4)(5)$.
Then $\sigma_2$ corresponds to the empty 0--1 tableau $\varphi=\varnothing$,
and we successively compute $S=[6]-Supp(53642)=\{ 1\}$,
$\Psi_{\{ 1\}}(w_aw_b)=1234566154$, and  $\tilde w_b=66154$.

Let $i_0$ denote
the length of the shortest column in $\varphi$, $1\leq i_0\leq cyc(\sigma_2)$. 
If $\varphi=\varnothing$, let $i_0=\infty$. Likewise, let $h_0$ denote the
smallest integer, $1\leq h_0\leq cyc(\sigma_2)$, 
such that $\tilde b_{h_0}<h_0$. If no
such $\tilde b_i$ exists, set $h_0=\infty$.

There are two cases: $i_0\geq h_0$ or $i_0< h_0$.
If $i_0\geq h_0$, then delete $\tilde b_{h_0}$ from the word $\tilde w_b$
and add a column of length $(h_0-1)$ to $\varphi$, with a 1 in position 
$\tilde b_{h_0}$, from bottom to top, thus obtaining a new pair 
$(\tilde w_b',\varphi')$. Since the letter removed from $\tilde w_b$ 
is at most equal to $(cyc(\sigma_2)-1)< (n_2-|B_2|)<(n_1+1)$,
$w_a\tilde w_b'$ is still an RG-function of maximum $(n_3+s)$, and the new
$i_0$ associated to $\varphi'$ is equal to $(h_0-1)$. 
$\Phi_5((B_i,\sigma_i);\pi)$ is then obtained by applying $\Psi_S^{-1}$ to
$w_a{\tilde w_b'}$ and by decoding the 0--1 tableau $\varphi'$.

In the above example, $i_0=\infty$ and $h_0=3$. Hence 
$\varphi'=\vcenter{\vbox{\btableaux{l} \plbox{0}\nbox\plbox{1}\nbox
\etableaux}}$\  (corresponding to the new permutation $\sigma_2=
(1)(2,3)(4)(5)$) and 
 $\tilde w_b'=6654$. From 
$\Psi_{\{ 1\}}^{-1}(w_a\tilde w_b)=123456165$, we get
$\Phi_5((B_i,\sigma_i);\pi)$ equals 
$B_1=B_2=B_3=\varnothing$, and 
$\pi=(\lo 1)(2,3)| (\lo 2)  (\hi 5)|(\lo 3)(\hi 2)|(\lo 4)(\hi 4)|
(\lo 5)(5)(\hi 1)|\break(1)(4)(\hi 3)$.

If $i_0<h_0$, the image of $((B_i,\sigma_i);\pi)$ is defined in the 
obvious way so that $\Phi_5$ is an involution. The proof that $\Phi_5$
is  weight-preserving and sign-reversing 
is quite straight-forward, and the details will
be left to the reader. It remains  to show
that $\Phi_5$ is well-defined. Remember that if $((B_i,\sigma_i);\pi)
\in Fix\Phi_4$, we must have $Supp(w_a)\subseteq Supp(w_bw_c)$. We 
have to show that $\Phi_5$ preserves this property. 
What complicates matters  is the application of $\Psi_S$ and
$\Psi_S^{-1}$ to the RG-functions $w_aw_b$. In Lemma 5, we explicitly
find the set of images $w_a\tilde w_b$ (which we will denote by 
$\tilde W(S)$) of all possible $w_aw_b$ under $\Psi_S$. We will
then show that the deletion or insertion of a letter whose value is 
strictly less than its position in $\tilde w_b$ yields  new
RG-functions $w_a\tilde w_b'$ which remain in the set $\tilde W(S)$.

Fix $n_3\geq n_1\geq n_2\geq 0$, $0\leq t\leq n_2$, and  
$0\leq s\leq n_1+n_2-n_3$. Let 
$S\subseteq [n_3+s]$ such that $|S|\leq s$, and  fix $w_a=12\ldots n_1$.
We denote by 
$$
\align
W(S)=\{&w_b\mid w_aw_b\in RG(n_1+n_2-t,n_3+s), \text{ and }\\
      &\, [n_1]\subseteq ([n_3+s]-S)\cup Supp(w_b)\},\\
\tilde W(S)=\{& \tilde w_b\mid \Psi_S(w_aw_b)=w_a\tilde w_b \text{ for }
               w_b\in W(S)\},\\
\endalign
$$
and
$$
\align
w_aW(S)=\{&w_aw_b\mid w_b\in W(S)\},\\
w_a\tilde W(S)=\{&w_a\tilde w_b\mid \tilde w_b\in\tilde W(S)\}.\\
\endalign
$$

In particular, when $|S|=s$, if $w_c$ is a word containing the letters
in $([n_3+s]-S)$ in any order, with no repetition, and $(B_2,\sigma_2)$ is
a partial permutation of $\{2\}\times[n_2]$ with $cyc(\sigma_2)=n_2-t$,
$W(S)$ contains all possible words $w_b$ such that $w=12\ldots n_1w_bw_c$
is the RG-function associated to some $((B_i,\sigma_i);\pi)\in Fix\Phi_4$
having these fixed $(B_2,\sigma_2)$ and $w_c$.

\proclaim{Lemma 5} (characterization of \, $\tilde W(S)$)\enspace
Let $S\subseteq [n_3+s]$ such that $|S|\leq s$. The set $\tilde W(S)$
depends only upon the cardinality $j=|S\cap [n_1]|$. More precisely, we have
\item{(i)} $\tilde W(S)=\tilde W(S\cap [n_1])$,
\item{(ii)} If $j=0$, $\tilde W(\varnothing)=W(\varnothing)$, and $\tilde w_b
\in \tilde W(\varnothing)$ has the following form:
$$
\tilde w_b={\undersetbrace \text{entries }\leq\,n_1\to{* \ldots *}}
(n_1+1) {\undersetbrace \leq\,(n_1+1)\to{* \ldots *}}
(n_1+2)\ldots (n_3+s-1){\undersetbrace \leq\,(n_3+s-1)\to{* \ldots *}}
(n_3+s){\undersetbrace \leq\,(n_3+s)\to{* \ldots *}}.
\tag 5.2
$$
\item{(iii)} If $j=1$, then $\tilde W(\{i\})=\tilde W(\{ 1\})$ is obtained
from $\tilde W(\varnothing)$ by keeping only the words $\tilde w_b$ of the 
form (5.2) such that one of the stars $*$ is set to its maximum and the
maximum value of all the stars to its right is lowered by 1. So
any $\tilde w_b$ has the form
$$
\align
\tilde w_b=&{\undersetbrace\text{entries } \leq\,n_1\to{* \ldots *}}
(n_1+1) {\undersetbrace \leq\,(n_1+1)\to{* \ldots *}}
(n_1+2)\ldots
(n_1+h){\undersetbrace \leq\,(n_1+h)\to{* \ldots *}}
(n_1+h){\undersetbrace \leq\,(n_1+h-1)\to{* \ldots *}}
\\ &(n_1+h+1)\ldots
(n_3+s-1){\undersetbrace \leq\,(n_3+s-2)\to{* \ldots *}}
(n_3+s){\undersetbrace \leq\,(n_3+s-1)\to{* \ldots *}}.\tag 5.3\\
\endalign
$$
\item{(iv)} If $j\geq 2$, then $\tilde W(S)=\tilde W(\{1,2,\ldots, j\})$ 
is obtained from $\tilde W(\{1,2,\ldots, j-1\})$ 
by the same construction as the one described in
(iii).
\endproclaim

\demo{Proof}\hfill\break
(i). First we show that $\tilde W(S)=\tilde W(S\cap [n_1])$. From the
definition of $W(S)$, it is clear that $W(S)=W(S\cap [n_1])$.
Moreover, if
$S=\{ s_1<\ldots <s_j<s_{j+1}<\ldots <s_n\}$, where $s_j\leq n_1$
and $s_{j+1}> n_1$, since $\Psi_{S\setminus [n_1]}$ is a bijection
on $RG(n_1+n_2-t, n_3+s)$, preserving first occurrences and leaving 
all letters $\leq n_1$ fixed, we must have
$$
\Psi_{S\setminus [n_1]}(w_aW(S))=w_aW(S).
$$
Therefore,
$$
\align
w_a\tilde W(S) &=\Psi_{S}(w_aW(S))
= \Psi_{S\cap [n_1]}\circ\Psi_{S\setminus [n_1]}(w_aW(S))\\
&=  \Psi_{S\cap [n_1]}(w_aW(S\cap[n_1]))=w_a\tilde W(S\cap[n_1]).\\
\endalign
$$

\noindent (ii). If $j=0$, $\Psi_{\varnothing}$ is the identity map and
$$
\tilde W(\varnothing)=W(\varnothing)=
\{w_b\mid 12\ldots n_1w_b\in RG(n_1+n_2-t,n_3+s)\},
$$
in which typical elements (tails of RG-functions) are given by (5.2).

\noindent (iii). If $j=1$, suppose $S=\{i\}$, $1\leq i\leq n_1$. Then 
$$
W(S)=\{ w_b\mid w_aw_b\in RG(n_1+n_2-t,n_3+s), 
\text{ and } i\in Supp(w_b)\}.
$$
Let $w_b\in W(S)$ and suppose the rightmost occurrence of $i$ lies in
position $p$ of $w_b$, between the first occurrence of $(n_1+h)$ and
the first occurrence of $(n_1+h+1)$. Thus $w_aw_b$ has the form
$$
\align
w_a w_b=& 1 2\ldots n_1\,{\undersetbrace\text{entries }
 \leq\,n_1\to{* \ldots *}}
(n_1+1) {\undersetbrace \leq\,(n_1+1)\to{* \ldots *}}
(n_1+2)\ldots
(n_1+h){\undersetbrace \leq\,(n_1+h)\to{* \ldots *}}
{\undersetbrace{\scriptstyle \text{position }\atop \scriptstyle (n_1+p)}\to i}
{\undersetbrace{\scriptstyle \leq\,(n_1+h),\atop
\scriptstyle \text{entries }
 \not= i}\to{* \ldots *}}
\\ &(n_1+h+1)\ldots
(n_3+s-1){\undersetbrace {\scriptstyle\leq\,(n_3+s-1),\atop\scriptstyle 
\not = i}\to{* \ldots *}}
(n_3+s){\undersetbrace {\scriptstyle\leq\,(n_3+s),\atop\scriptstyle 
\not = i}\to{* \ldots *}}.
\tag 5.4\\
\endalign
$$

Apply $\Psi_{\{i\}}=\Psi_{n_3+s}\circ\Psi_{n_3+s-1}\circ\ldots
\circ \Psi_{i}$ to $w_a w_b$. The last occurrence of $i$ in $w_aw_b$
(in position $(n_1+p)$) lies to the right of the first occurrence of
$(i+1)$ (case (i) in the definition of $\Psi_m$), 
so it is changed to $(i+1)$ by $\Psi_i$, and any $(i+1)$ to its
right is changed to $i$. Thus the last occurrence of $(i+1)$ in 
$\Psi_i(w_aw_b)$ now appears in position $(n_1+p)$, again to the right
of the first occurrence of $(i+2)$. So all $(i+2)$'s
to its right are changed to $(i+1)'s$ by $\Psi_{i+1}$, the $(i+1)$ in
position $(n_1+p)$ is switched to $(i+2)$, and every other letter remain
fixed.

The same argument applies until we reach $\Psi_{n_1+h}$. At this point in
$\Psi_{n_1+h-1}\circ\ldots\circ\Psi_i(w_aw_b)$, there is a $(n_1+h)$ in
position $(n_1+p)$ and no occurrence of $(n_1+h)$ to its right. 
This means that there are no letters $(n_1+h)$ to the right of the first
occurrence of $(n_1+h+1)$ (case (ii) in the definition of $\Psi_m$).
Hence $\Psi_{n_1+h}$ changes every occurrence
of $(n_1+h+1)$, except for the first one, to $(n_1+h)'s$, and fixes
everything else. Once again in the RG-function obtained, 
there are no occurrences of $(n_1+h+1)$
to the right of the first occurrence of $(n_1+h+2)$. It is clear
that by applying successively  $\Psi_{n_1+h+1},\ldots,\Psi_{n_3+s}$
respectively, we will get $\Psi_{\{i\}}(w_a w_b)$ exactly of the
form (5.3).
This shows that the set defined in (iii) is equal to $\tilde W(\{ i\})$.
Note that the definition of the set $\tilde W(\{ i\})$ is independent
of the actual value of $i$, so $\tilde W(\{ i\})=\tilde W(\{ 1\})$.

\noindent (iv). 
The proof is an easy induction based on the proof of (iii).
Note that if $S=\{s_1<s_2<\ldots<s_j\}, s_j\leq n_1$, the positions of the
last occurrences of  $s_1,s_2,\ldots,s_j$ respectively in $w_aw_b$
correspond exactly to the positions of the stars successively fixed
to their maximum in $\Psi_s(w_aw_b)$.\hfill $\square$
\enddemo

We can show now that $\Phi_5$ is well-defined.

Let $\tilde w_b\in\tilde W(\{1,2,\ldots, j\})$. The letters of $\tilde w_b$ can be
divided into two categories: the {\it fixed letters} (first occurrences
of $(n_1+1)$ up to $(n_3+s)$, and $j$ stars that were fixed to their
maximum in the construction described in the preceding lemma), and
the {\it free letters} (corresponding to stars in the description 
of $\tilde w_b$ in Lemma 5). So in order to be in 
$\tilde W(\{1,2,\ldots, j\})$,
a word  $\tilde w_b$ must have  $(n_3+s-n_1+j)$ fixed letters
(appearing in some fixed relative order), and possibly some free letters, 
depending on its length.

On one hand, note that the fixed letters of $\tilde w_b$ are always greater
or equal to their positions in $\tilde w_b$. Indeed, we have already
seen that the first occurrences were necessarily greater than their
position $p$ ($(n_3+s)\geq (n_1+1)>n_2\geq p$). As for the $j$
stars fixed to their maximum, the way to minimize their value in the
construction of Lemma 5 is to fix them successively by increasing order
of their positions. Then, if they all lie before the first
occurrences of $(n_1+1)$ up to $(n_3+s)$, the $j$-th star fixed will
have minimum value $(n_1-j+1)$, and the rightmost position 
where it can be located
is, for example, the one in the following word:
$$
\tilde w_b={\undersetbrace \text{entries}\leq\, n_1\to{**\ldots *}}\, n_1\,
(n_1-1)\ldots (n_1-j+1)(n_1+1)(n_1+2)\ldots (n_3+s).
$$
But from the relations $n_3\geq n_2$, $t\geq 0$,  and
$j\leq s$, we deduce that its position $p$,
$$
p=|\tilde w_b|-(n_3+s-n_1)=n_1+(n_2-n_3)-t-s\leq n_1-j+1.
$$
Therefore in that case, all fixed stars are greater or equal to their
positions. More generally, if a fixed star is rather located to the right of a 
first occurrence, its value is increased by one, 
so the letter  remains greater or equal to its position.

On the other hand, 
note  that the allowed maxima for the free letters are also greater
or equal to their positions in $\tilde w_b$. The same type of argument
(with same inequalities) applies. Details are left to the reader.

Now, the
``involutive step'' of  $\Phi_5$ was to add or to delete a letter
from $\tilde w_b$, and  this letter had the property of 
being strictly smaller
than its position in $\tilde w_b$. 

If the involutive step deleted  a letter from  $\tilde w_b$
($w_a\tilde w_b\in RG(n_1+n_2-t, n_3+s)$), then it had to be one of
its free letters because the fixed ones are greater or equal 
to their positions. Therefore the new  $\tilde w_b'$ obtained is
in the set $\tilde W(\{1,2,\ldots, j\})$ 
(with $w_a\tilde w_b'\in RG(n_1+n_2-(t+1), n_3+s)$).
Likewise, if the involutive step added a letter to $\tilde w_b$,
the new letter is in the right range to be considered  a free letter,
and the fixed letters (and their relative order) 
are not modified, so the new $\tilde w_b'$ is
in the set $\tilde W(\{1,2,\ldots, j\})$ as well (with $w_a\tilde w_b'\in 
RG(n_1+n_2-(t-1), n_3+s)$).

\subheading{Fixed  points for $\Phi_5$}
The fixed points of  $\Phi_5$ correspond to the case $i_0=h_0=\infty$.
Clearly, 
$$
\align
Fix\Phi_5=\{&((B_i,\sigma_i);\pi)\in Fix\Phi_4\mid \sigma_2=Id \text{ and
for } S=[\#blocks(\pi)]-Supp(w_c), \\ &\text{ the word }\tilde w_b \text{ in }
\Psi_S(w_aw_b)=w_a\tilde w_b
\text{ has its $i$-th letter }\geq i,\forall i\}.\\
\endalign
$$


\subheading{6. 
Combinatorial evaluation of $L_q(C_{n_1}(x)C_{n_2}(x)C_{n_3}(x))$}

An expression of  $L_q(C_{n_1}(x)C_{n_2}(x)C_{n_3}(x))$ can now be computed by
$\omega_q$-counting of the remaining fixed points $Fix\Phi_5$.
More precisely, $((B_i,\sigma_i);\pi)\in Fix\Phi_5$ if and only if

\noindent{Fix.1}\enspace $B_1=B_3=\varnothing$,

\noindent{Fix.2}\enspace $\sigma_i=Id$ for $i=1,2,3$,

\noindent{Fix.3}\enspace $w_a$ (respectively $w_c$) has all distinct letters 
and $Supp(w_a)\subseteq Supp(w_bw_c)$ (respec- \break
\phantom{\noindent{Fix.3}\enspace}tively $Supp(w_c)\subseteq Supp(w_aw_b)$),

\noindent{Fix.4}\enspace for $S=[\#blocks(\pi)]-Supp(w_c)$, the word
$\tilde w_b=\tilde b_1\tilde b_2\ldots\tilde b_{n_2-|B_2|}$ in
$\Psi_S(w_aw_b)=$  \break
\phantom{\noindent{Fix.4}\enspace}$w_a\tilde w_b$ has all $\tilde b_i$'s $\geq i$,
where $\Psi_S$ was defined in Lemma 4.

Clearly, for such elements, the weight (as was defined in (4.8)) reduces to
$$
\omega_q((B_i,\sigma_i):\pi)=(-1)^{|B_2|}q^{inv(B_2)+rs(\pi)}
a^{|B_2|+\#blocks(\pi)}.
\tag 6.1
$$

By $\omega_q$-counting this fixed point set, we will show that
$$
\align
&L_q(C_{n_1}(x)C_{n_2}(x)C_{n_3}(x))=
\sum_{l=0}^{n_1+n_2-n_3}
\sum_{s=0}^l\sum_{j=0}^sa^{n_3+l}(-1)^{l-s}q^{L}
[n_3]!_q
\left[\matrix n_2\\l-s\endmatrix\right]_q\tag 6.2\\ &
 \left[\matrix n_1\\j\endmatrix\right]_q\left[\matrix n_3-n_1+s\\
s-j\endmatrix\right]_q
\left[\matrix n_2-l+s\\n_3-n_1+s\endmatrix\right]_q
{[j]!_q[n_1-j]!_q\over [n_3-n_2+l]!_q}
\left[\matrix n_1+n_2-n_3-l\\j\endmatrix\right]_q,\\
\endalign
$$
where
$$
\align
L=&\binom{n_1}{2}+\binom{l-s}{2}+
j(-n_3-s+1)+\binom{j}{2}-\binom{s-j}{2}-(s-j)(n_3-n_1+j)\\&
+\binom{n_2-l+s}{2}+(n_3-n_1+s)(n_3-n_2+l)+j(n_3-n_2+l).\\
\endalign
$$
Evaluating  the $s$-sum by the $q$-binomial theorem (which has a
simple bijective proof) gives the right-hand side of (4.5) and thus Theorem 3.

The main difficulty here is to transpose the condition Fix.4 into
the $\omega_q$-counting. Using Lemmas 4 and 5, we will see that 
this corresponds  to the $q$-counting 
of some special sets
of RG-functions according to the statistic $ls$, which is the object of
Lemma 6.

Let us first group the elements of $Fix\Phi_5$ 
by powers of $a$. The power of $a$ ranges
from a minimum of $n_3$ (expressing the fact that $w_c$ has $n_3$ distinct
letters) to a maximum of $(n_1+n_2)$ (being the maximum value of 
$\{\text{max}(Supp(
w_aw_b))+|B_2|\}$). Now,
$$
\align
L_q(C_{n_1}C_{n_2}C_{n_3})
&=\sum_{l=0}^{n_1+n_2-n_3}a^{n_3+l}\sum_{s=0}^l(-1)^{l-s}q^{\binom{l-s}{2}}
 \left[\matrix n_2\\l-s\endmatrix\right]_q
\sum_{\pi} q^{rs(\pi)},\tag 6.3\\
\endalign
$$
where the last sum ranges over all partitions $\pi$ corresponding to
$((B_i,\sigma_i);\pi)\in Fix\Phi_5$ such that $\#blocks(\pi)=(n_3+s)$
and $B_2$ is any fixed subset of $\{2\}\times[n_2]$ of cardinality $(l-s)$.
The $s$-sum is the generating function for the subsets $B_2$, as was
established in \S2. But
$$
\align
rs(\pi) &= rs(w)
        = rs(w_c)+rs(w)|_{w_aw_b}\\
      &= rs(w_c)+ls(w_a\tilde w_b)-\sum_{u\in \left([n_3+s]-Supp(w_c)\right)}
     (n_3+s-u).\tag 6.4\\
\endalign
$$
Note that for any fixed set $Supp(w_c)$, there are no constraints on the
positions of the letters in $w_c$, so $rs(w_c)$ is simply the number of
inversions of the word $w_c$, whose distribution is mahonian (i.e.
the generating function equals $[n_3]!_q$).
From Lemma 5, we also know that the possible choices for $\tilde w_b$
only depend on the cardinality $j$ of the set 
$\left([n_3+s]-Supp(w_c)\right)\cap[n_1]$, not on the actual
set $Supp(w_c)$ itself. Hence, if we let 
$$
Fix\tilde W(j)=\{ \tilde w_b\mid \tilde w_b=\tilde b_1\ldots \tilde b_{
n_2-l+s}\in \tilde W(\{1,2,\ldots, j\}) 
\text{ and } \tilde b_i\geq i,\forall i\},
$$
where $\tilde W(\{1,2,\ldots, j\})$ was characterized in  Lemma 5, we get
that the last sum on the right-hand side of (6.3) equals
$$
\align
\sum_{\pi} q^{rs(\pi)}
&=
[n_3]!_q \, 
\sum_{j=0}^s\, \sum_{\tilde w_b\in Fix\tilde W(j)}q^{ls(w_a\tilde w_b)}\\ &
q^{j(-n_3-s+1)+\binom{j}{2}}
\left[\matrix n_1\\j\endmatrix\right]_q
q^{-\binom{s-j}{2}-(s-j)(n_3-n_1+j)}
 \left[\matrix n_3-n_1+s\\s-j\endmatrix\right]_q.
\tag 6.5\\
\endalign
$$
Finally, we show that
\proclaim{Lemma 6} If $w_a=12\ldots n_1$, $\tilde w_b=\tilde b_1\ldots
\tilde b_{n_2-l+s}$ and $max\left(Supp(w_a\tilde w_b)\right)=n_3+s$, then
$$
\sum_{\tilde w_b\in Fix\tilde W(j)}q^{ls(w_a\tilde w_b)}=
q^{A}
\left[\matrix n_2-l+s\\n_3-n_1+s\endmatrix\right]_q
{[j]!_q[n_1-j]!_q\over [n_3-n_2+l]!_q}
\left[\matrix n_1+n_2-n_3-l\\j\endmatrix\right]_q,\tag 6.6
$$
where
$$
A=\binom{n_1}{2}+\binom{n_2-l+s}{2}+(n_3-n_1+s+j)(n_3-n_2+l).
$$
\endproclaim

\demo{Proof} Note that the statistic $ls$ of any RG-function is just the 
sum of the values of the letters minus one, so
$$
ls(w_a\tilde w_b)=\binom{n_1}{2}+\sum_{i=1}^{n_2-l+s}(\tilde b_i-1).
$$

To visualize more easily where the various factors of (6.6) come from, let
us encode $\tilde w_b$ as a 0--1 tableau $\varphi$ in the following manner:
start with a $(n_3+s)\times(n_2-l+s)$ rectangular Ferrers diagram. Fill it
with a 1 in position $j$ (from bottom to top) of column $i$ if $\tilde b_i=j$,
and with 0's elsewhere. For example, if $(n_3+s)=8$, $n_1=6=(n_2-l+s)$ and
$\tilde w_b=175787$, $\varphi$ is the  0--1 tableau on the left of Figure 3.

\midinsert
$$
\varphi\,=\,
\vcenter{
\vbox{\btableaux{l}
\plbox{0}\plbox{0}\plbox{0}\plbox{0}\plbox{1}\plbox{0}\nbox
\plbox{0}\plbox{1}\plbox{0}\plbox{1}\plbox{0}\plbox{1}\nbox
\plbox{0}\plbox{0}\plbox{0}\plbox{0}\plbox{0}\plbox{0}\nbox
\plbox{0}\plbox{0}\plbox{1}\plbox{0}\plbox{0}\shbox{0}\nbox
\plbox{0}\plbox{0}\plbox{0}\plbox{0}\shbox{0}\shbox{0}\nbox
\plbox{0}\plbox{0}\plbox{0}\shbox{0}\shbox{0}\shbox{0}\nbox
\plbox{0}\plbox{0}\shbox{0}\shbox{0}\shbox{0}\shbox{0}\nbox
\plbox{1}\shbox{0}\shbox{0}\shbox{0}\shbox{0}\shbox{0}\nbox
\etableaux}}
\quad\longrightarrow\quad
\tilde\varphi_{typ}\,=\,
\vcenter{
\vbox{\btableaux{l}
\plbox{0}\plbox{0}\plbox{0}\plbox{0}\plbox{1}\plbox{{$*$}}\nbox
\plbox{0}\plbox{1}\plbox{{$*$}}\plbox{{$*$}}\shbox{0}\plbox{{$*$}}\nbox
\plbox{{$*$}}\shbox{0}\plbox{{$*$}}\plbox{{$*$}}\shbox{0}\plbox{{$*$}}\nbox
\plbox{{$*$}}\shbox{0}\plbox{{$*$}}\plbox{{$*$}}\shbox{0}\nbox
\plbox{{$*$}}\shbox{0}\plbox{{$*$}}\plbox{{$*$}}\nbox
\plbox{{$*$}}\shbox{0}\plbox{{$*$}}\nbox
\plbox{{$*$}}\shbox{0}\nbox
\plbox{{$*$}}\nbox
\etableaux}}
$$
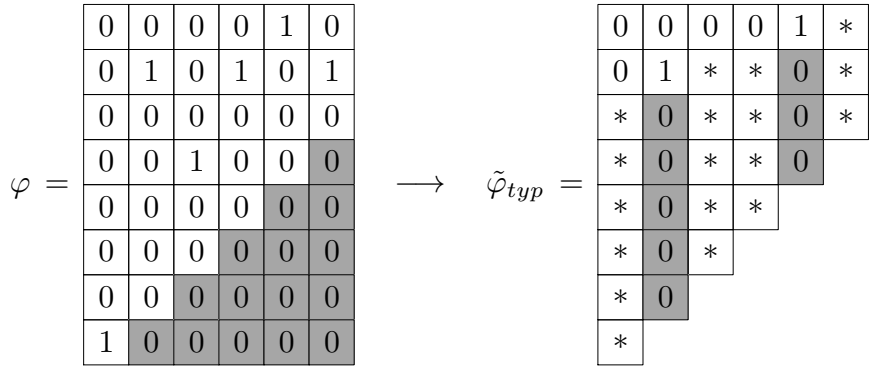
\centerline{Figure 3: Encoding of $\tilde w_b$ as a 0--1 tableau}
\endinsert

Obviously, we have $\sum(\tilde b_i-1)=inv(\varphi)$. Note also that the
0's in the shaded staircase shape of $\varphi$ in Figure 3 always count
as inversions, expressing the fact that $\tilde b_i\geq i$. They account
for the factor $q^{\binom{n_2-l+s}{2}}$ in (6.6). We can therefore drop
them from $\varphi$ without loss of generality, and compute the inversion
number of the reduced 0--1 tableau $\tilde\varphi$. We now use Lemma 5 to
characterize the possible fillings of $\tilde\varphi$ according to $j$.
\enddemo 

\subheading{Case 1: $j=0$}
From Lemma 5 (ii), $\tilde w_b\in\tilde W(\varnothing)$ simply means that it
is the tail of an RG-function. Thus the only restrictions on $\tilde w_b$
are that the first occurrences of $(n_1+1), (n_1+2),\ldots,(n_3+s)$
appear in the right order.

If we set $x=(n_3+s)$, $y=n_1$, and $z=(n_2-l+s)$, in the context of
0--1 tableaux, we want to $q$-count all 0--1 tableaux $\tilde \varphi$
with $z$ columns of lengths $x,(x-1),\ldots, (x-z+1)$ respectively, such that
when we look at the top $(x-y)$ rows of  $\tilde \varphi$ from left to
right, the leftmost 1 in any row must always occur before the ones in the
rows above it. Grouping the tableaux according to these leftmost occurrences
of 1's, we get ``typical'' 0--1 tableaux  $\tilde \varphi_{typ}$, 
corresponding exactly to the typical words $\tilde w_b$ described in 
(5.2) of Lemma 5. For instance, the typical 0--1 tableau 
$\tilde \varphi_{typ}$ containing our previous example is illustrated in
Figure 3 (stars $*$ correspond to possible positions of 1's).
Carrying out the $q$-counting, observe that
\item{1.1} Each column containing a number $m$ of stars contributes a 
factor $[m]_q$ to the $q$-counting of inversions. No matter which $(x-y)$
columns are chosen to be first occurrences of upper 1's, the number of
stars in the remaining columns is $y, (y-1), \ldots$ and $(x-z+1)$
respectively, contributing to an overall factor of 
$${[y]!_q\over [x-z]!_q}={[n_1]!_q\over [n_3-n_2+l]!_q}.$$
\item{1.2} The 0's below the leftmost occurrences of upper 1's (shaded in
Figure 3) form a partition $\mu$ with $(x-y)$ parts of length at least
$(x-z)$ and at most $y$, determined by the positions of the first
occurrences. Summing over all possible choices, it contributes a factor
$$q^{(x-y)(x-z)}
\left[\matrix z\\ x-y\endmatrix\right]_q
=q^{(n_3-n_1+s)(n_3-n_2+l)}
\left[\matrix n_2-l+s\\ n_3-n_1+s\endmatrix\right]_q.$$
\subheading{Case 2: $j\geq 1$}
Recall that Lemma 5 (iii) and (iv) provides a method to construct all the
elements of $\tilde W(\{1,2,\ldots, j\})$ uniquely from 
$\tilde W(\varnothing)$. In the 
0--1 tableau context, if we extract only the cells filled with stars in
$\tilde\varphi_{typ}$ (hence obtaining a tableau $\tilde\psi_{typ}$
with $(n_1+n_2-n_3-l)$ columns of lengths $n_1,(n_1-1),\ldots, 
(n_3-n_2+l+1)$ respectively), the manipulation described in Lemma 5 (iii)
corresponds to replacing the top star of a column   by a 1, and all the
top stars to its right and the stars below it by a 0. Repeating this
procedure $j$ times and reinserting the columns of $\tilde\psi_{typ}$
in $\tilde\varphi_{typ}$ yields to ``typical'' 0--1 tableaux that correspond
to the elements of $Fix\tilde W(j)$. For example, Figure 4 shows the
above manipulations on the third and the first columns respectively
of the stars extracted from
$\tilde\varphi_{typ}$ of Figure 3.

\midinsert
$$
\vcenter{
\vbox{\btableaux{l}
\plbox{{$*$}}\plbox{{$*$}}\plbox{{$*$}}\plbox{{$*$}}\nbox
\plbox{{$*$}}\plbox{{$*$}}\plbox{{$*$}}\plbox{{$*$}}\nbox
\plbox{{$*$}}\plbox{{$*$}}\plbox{{$*$}}\plbox{{$*$}}\nbox
\plbox{{$*$}}\plbox{{$*$}}\plbox{{$*$}}\nbox
\plbox{{$*$}}\plbox{{$*$}}\nbox
\plbox{{$*$}}\nbox
\etableaux}}
\quad\longrightarrow\quad
\vcenter{
\vbox{\btableaux{l}
\plbox{{$*$}}\plbox{{$*$}}\plbox{{1}}\plbox{{0}}\nbox
\plbox{{$*$}}\plbox{{$*$}}\plbox{{0}}\plbox{{$*$}}\nbox
\plbox{{$*$}}\plbox{{$*$}}\plbox{{0}}\plbox{{$*$}}\nbox
\plbox{{$*$}}\plbox{{$*$}}\plbox{{0}}\nbox
\plbox{{$*$}}\plbox{{$*$}}\nbox
\plbox{{$*$}}\nbox
\etableaux}}
\quad\longrightarrow\quad
\vcenter{
\vbox{\btableaux{l}
\plbox{{1}}\plbox{{0}}\plbox{{1}}\plbox{{0}}\nbox
\plbox{{0}}\plbox{{$*$}}\plbox{{0}}\plbox{{0}}\nbox
\plbox{{0}}\plbox{{$*$}}\plbox{{0}}\plbox{{$*$}}\nbox
\plbox{{0}}\plbox{{$*$}}\plbox{{0}}\nbox
\plbox{{0}}\plbox{{$*$}}\nbox
\plbox{{0}}\nbox
\etableaux}}
\quad\Longrightarrow\quad
\vcenter{
\vbox{\btableaux{l}
\plbox{0}\plbox{0}\plbox{0}\plbox{0}\plbox{1}\shbox{{0}}\nbox
\plbox{0}\plbox{1}\shbox{{0}}\shbox{{1}}\plbox{0}\shbox{{0}}\nbox
\shbox{{1}}\plbox{0}\shbox{{$*$}}\shbox{{0}}\plbox{0}\shbox{{$*$}}\nbox
\shbox{{0}}\plbox{0}\shbox{{$*$}}\shbox{{0}}\plbox{0}\nbox
\shbox{{0}}\plbox{0}\shbox{{$*$}}\shbox{{0}}\nbox
\shbox{{0}}\plbox{0}\shbox{{$*$}}\nbox
\shbox{{0}}\plbox{0}\nbox
\shbox{{0}}\nbox
\etableaux}}
$$
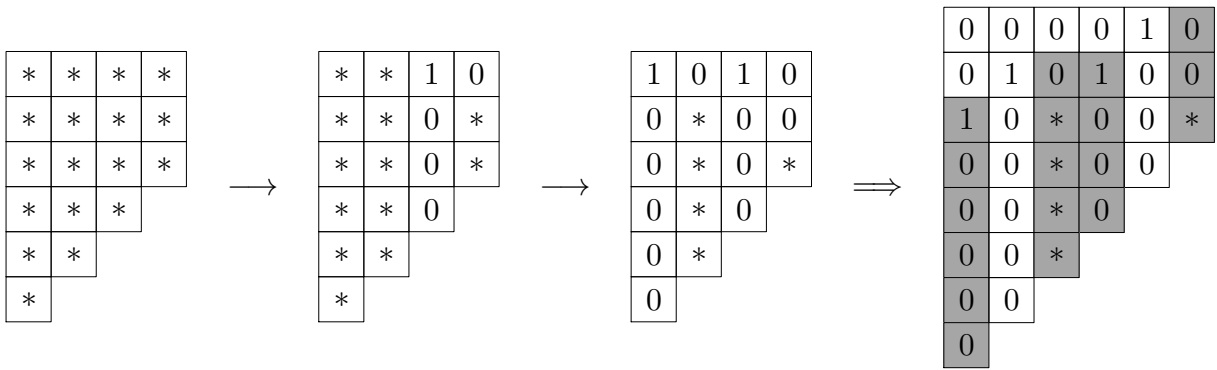
\centerline{Figure 4: Manipulations of Lemma 5 in the context of 0--1 tableaux}
\endinsert

Proceeding to $q$-counting, the part (1.2) of case 1 is left unchanged and
the part (1.1) is replaced by the contribution  of the different choices
of $\tilde\psi_{typ}$. But observe that
\item{2.1} All the 0--1 tableaux  $\tilde\psi_{typ}$ such that a star has
been changed to a 1 in columns $c_1,c_2,\ldots$ and $c_j$ contribute to
$[j]!_q$ times the $q$-counting of the 0--1 tableaux  $\tilde\psi_{typ}$
such that this procedure was done in increasing order of the $c_i$'s.
Therefore, we can restrict to this latter case. This explains the factor
$[j]!_q$ in (6.6).
\item{2.2} It is not hard to see that in that case, we are $q$-counting all
0--1 tableaux  $\tilde\psi$ containing $(n_1+n_2-n_3-l)$ columns of lengths
$n_1,(n_1-1), \ldots, (n_3-n_2+l+1)$ respectively, such that when we look
at the top $j$ rows from right to left, the rightmost 1 in any row has
to occur before the ones in the rows above it. There is a simple
weight-preserving bijection between this class of 0--1 tableaux and the
one that was $q$-counted in case 1, for $x=n_1$, $y=(n_1-j)$ and
$z=(n_1+n_2-n_3-l)$ (this class is defined by interchanging ``left'' and
 ``right''). Given $\tilde\psi$ in the first
class of 0--1 tableaux, just leave all the 1's below the $j$-th row
fixed and ``reverse the order'' of the 1's in the top $j$ rows, within
the columns where they appear. Figure 5
gives an example of this for $j=2$.

\midinsert
$$
\vcenter{
\vbox{\btableaux{l}
\plbox{0}\plbox{1}\plbox{1}\plbox{0}\plbox{0}\plbox{0}\nbox
\plbox{0}\plbox{0}\plbox{0}\plbox{0}\plbox{0}\plbox{1}\nbox
\shbox{0}\shbox{0}\shbox{0}\shbox{1}\shbox{0}\shbox{0}\nbox
\shbox{0}\shbox{0}\shbox{0}\shbox{0}\shbox{1}\nbox
\shbox{0}\shbox{0}\shbox{0}\shbox{0}\nbox
\shbox{1}\shbox{0}\shbox{0}\nbox
\shbox{0}\shbox{0}\nbox
\shbox{0}\nbox
\etableaux}}
\quad\longleftrightarrow\quad
\vcenter{
\vbox{\btableaux{l}
\plbox{0}\plbox{0}\plbox{1}\plbox{0}\plbox{0}\plbox{1}\nbox
\plbox{0}\plbox{1}\plbox{0}\plbox{0}\plbox{0}\plbox{0}\nbox
\shbox{0}\shbox{0}\shbox{0}\shbox{1}\shbox{0}\shbox{0}\nbox
\shbox{0}\shbox{0}\shbox{0}\shbox{0}\shbox{1}\nbox
\shbox{0}\shbox{0}\shbox{0}\shbox{0}\nbox
\shbox{1}\shbox{0}\shbox{0}\nbox
\shbox{0}\shbox{0}\nbox
\shbox{0}\nbox
\etableaux}}
$$
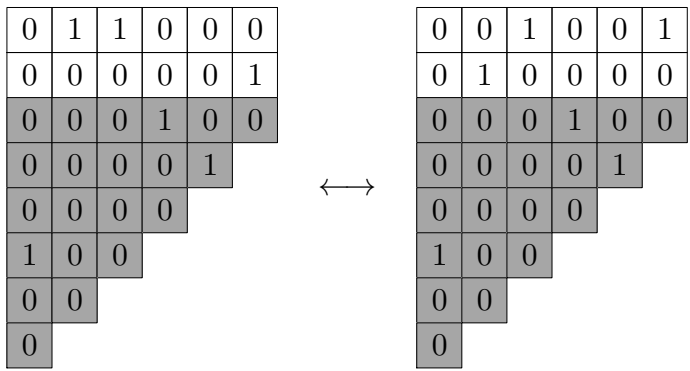
\centerline{Figure 5: Bijection between two classes of 0--1 tableaux}
\endinsert

This is clearly an involution that preserves the number of 0's below 1's.
Therefore, we can simply use case 1 to compute the $q$-contribution of the
$\tilde\psi$'s. We obtain
\TagsOnRight
$$
q^{j(n_3-n_2+l)}{[n_1-j]!_q\over[n_3-n_2+l]!_q}
\left[\matrix n_1+n_2-n_3-l\\j\endmatrix\right]_q.\tag"$\square$"
$$
\TagsOnLeft

Finally, putting together Lemma 6, 
identities (6.5) and (6.3) yields identity (6.2),
thus completing the proof of Theorem 3. 

Note that if we take $n_2=0$ and apply $\Phi_1$ and $\Phi_2$ to 
$L_q(n_1,0,n_3)$ (assuming $n_3\geq n_1$), the set $Fix\Phi_2$ is easily seen
to be empty unless $n_1=n_3$, in which case it can be proven to be
$\omega_q$-counted by the right-hand side of (4.2) ($n=n_1, m=n_3$),
thus proving orthogonality and Theorem 2.
These results can also be obtained by applying directly $\Phi_3$ and 
$\Phi_4$ to the set $L_q(n_1,0,n_3)$, assuming this time $n_1\geq n_3$.
We also have a weight-preserving sign-reversing involution proving
orthogonality when the $q$-statistic for the moments is taken to be 
$lb$ instead of $rs$, but we do not know how to generalize it to the
linearization problem.

\proclaim{Corollary 7}
Let $n_1\geq n_2\geq \ldots \geq n_k$. The coefficient of the lowest
power of $a$, $a^{n_1}$ in $L_q(C_{n_1}C_{n_2}\ldots C_{n_k})$ is a
polynomial in $q$ with positive coefficients.
\endproclaim

\demo{Proof}
The proof of Theorem 3 can be generalized to a product of $k$ $q$-Charlier
polynomials, any additional color being treated as was color 2, the middle
color. It is
easy to see then that the fixed points contributing to the lowest power
of $a$ must have all $B_i=\varnothing$, and therefore have all
positive weights.\hfill$\square$
\enddemo

\proclaim{Corollary 8}
Let $n_3\geq n_1\geq n_2$. The coefficient of $a^{n_1+n_2-i}$ in
$L_q(C_{n_1}C_{n_2}C_{n_3})$ is equal to $(q-1)^{n_1+n_2-n_3-2i}$
times the coefficient of $a^{n_3+i}$, for
$0\leq i\leq\lfloor (n_1+n_2-n_3)/2\rfloor$.
\endproclaim

Our proof of Corollary 8 is analytical, but we would like to have
a combinatorial explanation of this ``symmetry'' property.

Note that $Fix\Phi_5$ is not an optimal set of fixed points,
in the sense that there are
still some terms that cancel each other when we proceed to
$\omega_q$-counting of $Fix\Phi_5$.
For example, for $n_1=n_2=n_3=2$, the two elements
of $Fix\Phi_5$ such that $B_2=\{2\}$, $w=12121$ and
$B_2=\varnothing$, $w=123123$ have weight $-a^3q^3$ and $a^3q^3$
respectively. However, we do not believe that an attempt to reduce
$Fix\Phi_5$ would be worthwhile.

\proclaim{Corollary 9}
Let $n_1\geq n_2\geq \ldots \geq n_k$. If
$q=1+r$, $L_q(C_{n_1}C_{n_2}\ldots C_{n_k})$ is a polynomial in $r$
with positive coefficients.
\endproclaim


\subheading{7. The classical $q$-Charlier polynomials}

We contrast the results of the previous sections with those for the 
classical $q$-Charlier polynomials \cite{11, p.187}
$$
c_n(x;a;q)= \ _2\phi_1(q^{-n},x;0;q,-q^{n+1}/a).
\tag7.1 
$$

The monic form of these polynomials, $cc_n(x;a;q)$ satisfies
$$
cc_{n+1}(x;a;q)=(x-b_n)cc_n(x;a;q)-\lambda_n cc_{n-1}(x;a;q),
$$
where
$$
\lambda_n=-aq^{1-2n}(1-q^{-n})(1+aq^{-n}),\quad b_n=aq^{-1-2n}+q^{-n}+aq^{-2n}-aq^{-n}.
$$
A calculation (see \cite{11, p.187}) shows that
 the moments for these polynomials are
$$
\mu_n=\prod_{i=1}^n (1+aq^{-i}).
$$

We need to rescale $x$ and $a$ so that $b_n$ and 
$\lambda_n$ are $q$-analogues of $a+n$ and $an$ respectively. If we put 
$x=1+z(1-q)$, and multiply $a$ by $(1-q)$, and call the resulting 
monic polynomials
$\hat C_n(z;a;q)$, the explicit formula from (7.1) is
$$
\hat C_n(z;a;q)=q^{-n^2}\sum_{k=0}^n\left[\matrix n\\k\endmatrix\right]_q
(-a)^{n-k}q^{\binom{k+1}2}
\prod_{i=0}^{k-1}(q^iz-[i]_q)
\tag7.2
$$
The three term recurrence relation coefficients are
$$
b_n=q^{-n}[n]_q(1+a(1-q)q^{-n})+aq^{-1-2n},\quad
\lambda_n=aq^{1-3n}[n]_q(1+a(1-q)q^{-n}).
\tag 7.3
$$
A calculation using the measure in \cite{11, p.187} gives
$$
\mu_n=\sum_{j=1}^n 
q^{-\binom{j}{2}-n}S_{1/q}(n,j)a^j.
\tag7.4
$$
Again we find $q$-Stirling numbers for the moments.
Zeng \cite{24} has also derived (7.2) and (7.3) from the continued
fraction for the moment  generating function.

We see that the individual terms in (7.3) do not have 
constant sign. This means that the Viennot theory must involve a 
sign-reversing involution for its combinatorial versions of (7.3)
and (7.4). Nonetheless we can give combinatorial interpretations of
(7.2) and (7.4), but have no perfect analog of Theorem 3.

\subheading{Acknowledgement} Theorems 1 and 2 were originally found
in joint work with Mourad Ismail.


\enddocument

%% file: texdraw.tex





 
\def\setRevDate $#1 #2 #3${\def\TeXdrawId{TeXdraw V1R4a <#2>}}
\setRevDate $Date: 1992/07/07 14:46:28 $

\chardef\catamp=\the\catcode`\@
\catcode`\@=11

\long                              
\def\centertexdraw #1{\hbox to \hsize{\hss
                                      \btexdraw #1\etexdraw
                                      \hss}}


\def\btexdraw {\x@pix=0             \y@pix=0
               \x@segoffpix=\x@pix  \y@segoffpix=\y@pix
               \t@exdrawdef
               \setbox\t@xdbox=\vbox\bgroup\offinterlineskip
                   \global\d@bs=0           
                   \t@extonlytrue           
                   \p@osinitfalse
                   \s@avemove \x@pix \y@pix 
                   \m@pendingfalse
                   \p@osinitfalse           
                   \p@athfalse}


\def\etexdraw {\ift@extonly \else
                 \t@drclose      
               \fi
               \egroup           
               \ifdim \wd\t@xdbox>0pt
                 \errmessage{TeXdraw box non-zero size,
                             possible extraneous text}%
               \fi
               \maxhvpos         
               \pixtodim \xminpix \l@lxpos  \pixtodim \yminpix \l@lypos
               \vbox {\vskip \vdrawsize
                      \t@xdinclude              
                      \vskip \l@lypos
                      \hbox {\hskip -\l@lxpos
                             \box\t@xdbox       
                             \hskip \hdrawsize
                             \hskip \l@lxpos}%
                      \vskip -\l@lypos\relax}}

%

\def\drawdim #1 {\def\d@dim{#1\relax}}


\def\setunitscale #1 {\edef\u@nitsc{#1}%
                      \realmult \u@nitsc  \s@egsc \d@sc}
\def\relunitscale #1 {\realmult {#1}\u@nitsc \u@nitsc
                      \realmult \u@nitsc \s@egsc \d@sc}
\def\setsegscale #1 {\edef\s@egsc {#1}%
                     \realmult \u@nitsc \s@egsc \d@sc}
\def\relsegscale #1 {\realmult {#1}\s@egsc \s@egsc
                     \realmult \u@nitsc \s@egsc \d@sc}

\def\bsegment {\ifp@ath
                 \f@lushbs
                 \f@lushmove
               \fi
               \begingroup
               \x@segoffpix=\x@pix
               \y@segoffpix=\y@pix
               \setsegscale 1
               \global\advance \d@bs by 1\relax}
\def\esegment {\endgroup
               \ifnum \d@bs=0
                 \writetx {es}%
               \else
                 \global\advance \d@bs by -1
               \fi}

\def\savecurrpos (#1 #2){\getsympos (#1 #2)\a@rgx\a@rgy
                         \s@etcsn \a@rgx {\the\x@pix}%
                         \s@etcsn \a@rgy {\the\y@pix}}%
\def\savepos (#1 #2)(#3 #4){\getpos (#1 #2)\a@rgx\a@rgy
                            \coordtopix \a@rgx \t@pixa
                            \advance \t@pixa by \x@segoffpix
                            \coordtopix \a@rgy \t@pixb
                            \advance \t@pixb by \y@segoffpix
                            \getsympos (#3 #4)\a@rgx\a@rgy
                            \s@etcsn \a@rgx {\the\t@pixa}%
                            \s@etcsn \a@rgy {\the\t@pixb}}

\def\linewd #1 {\coordtopix {#1}\t@pixa
                \f@lushbs
                \writetx {\the\t@pixa\space sl}}
\def\setgray #1 {\f@lushbs
                 \writetx {#1 sg}}
\def\lpatt (#1){\listtopix (#1)\p@ixlist
                \f@lushbs
                \writetx {[\p@ixlist] sd}}

\def\lvec (#1 #2){\getpos (#1 #2)\a@rgx\a@rgy
                  \s@etpospix \a@rgx \a@rgy
                  \writeps {\the\x@pix\space \the\y@pix\space lv}}
\def\rlvec (#1 #2){\getpos (#1 #2)\a@rgx\a@rgy
                   \r@elpospix \a@rgx \a@rgy
                   \writeps {\the\x@pix\space \the\y@pix\space lv}}
\def\move (#1 #2){\getpos (#1 #2)\a@rgx\a@rgy
                  \s@etpospix \a@rgx \a@rgy
                  \s@avemove \x@pix \y@pix}
\def\rmove (#1 #2){\getpos (#1 #2)\a@rgx\a@rgy
                   \r@elpospix \a@rgx \a@rgy
                   \s@avemove \x@pix \y@pix}

\def\lcir r:#1 {\coordtopix {#1}\t@pixa
                \writeps {\the\t@pixa\space cr}%
                \r@elupd \t@pixa \t@pixa
                \r@elupd {-\t@pixa}{-\t@pixa}}
\def\fcir f:#1 r:#2 {\coordtopix {#2}\t@pixa
                     \writeps {#1 \the\t@pixa\space fc}%
                     \r@elupd \t@pixa \t@pixa
                     \r@elupd {-\t@pixa}{-\t@pixa}}
\def\lellip rx:#1 ry:#2 {\coordtopix {#1}\t@pixa
                         \coordtopix {#2}\t@pixb
                         \writeps {\the\t@pixa\space \the\t@pixb\space el}%
                         \r@elupd \t@pixa \t@pixb
                         \r@elupd {-\t@pixa}{-\t@pixb}}
\def\larc r:#1 sd:#2 ed:#3 {\coordtopix {#1}\t@pixa
                            \writeps {\the\t@pixa\space #2 #3 ar}}


\def\ifill f:#1 {\writeps {#1 fl}}     
\def\lfill f:#1 {\writeps {#1 fp}}     



\def\htext #1{\def\testit {#1}%
              \ifx \testit\l@paren
                \let\next=\h@move
              \else
                \let\next=\h@text
              \fi
              \next {#1}}

\def\rtext td:#1 #2{\def\testit {#2}%
                    \ifx \testit\l@paren
                      \let\next=\r@move
                    \else
                      \let\next=\r@text
                    \fi
                    \next td:#1 {#2}}


\def\textref h:#1 v:#2 {\ifx #1R%
                          \edef\l@stuff {\hss}\edef\r@stuff {}%
                        \else
                          \ifx #1C%
                            \edef\l@stuff {\hss}\edef\r@stuff {\hss}%
                          \else  
                            \edef\l@stuff {}\edef\r@stuff {\hss}%
                          \fi
                        \fi
                        \ifx #2T%
                          \edef\t@stuff {}\edef\b@stuff {\vss}%
                        \else
                          \ifx #2C%
                            \edef\t@stuff {\vss}\edef\b@stuff {\vss}%
                          \else  
                            \edef\t@stuff {\vss}\edef\b@stuff {}%
                          \fi
                        \fi}

\def\avec (#1 #2){\getpos (#1 #2)\a@rgx\a@rgy
                  \s@etpospix \a@rgx \a@rgy
                  \writeps {\the\x@pix\space \the\y@pix\space (\a@type)
                            \the\a@lenpix\space \the\a@widpix\space av}}

\def\ravec (#1 #2){\getpos (#1 #2)\a@rgx\a@rgy
                   \r@elpospix \a@rgx \a@rgy
                   \writeps {\the\x@pix\space \the\y@pix\space (\a@type)
                             \the\a@lenpix\space \the\a@widpix\space av}}

\def\arrowheadsize l:#1 w:#2 {\coordtopix{#1}\a@lenpix
                              \coordtopix{#2}\a@widpix}
\def\arrowheadtype t:#1 {\edef\a@type{#1}}

\def\clvec (#1 #2)(#3 #4)(#5 #6)%
           {\getpos (#1 #2)\a@rgx\a@rgy
            \coordtopix \a@rgx\t@pixa
            \advance \t@pixa by \x@segoffpix
            \coordtopix \a@rgy\t@pixb
            \advance \t@pixb by \y@segoffpix
            \getpos (#3 #4)\a@rgx\a@rgy
            \coordtopix \a@rgx\t@pixc
            \advance \t@pixc by \x@segoffpix
            \coordtopix \a@rgy\t@pixd
            \advance \t@pixd by \y@segoffpix
            \getpos (#5 #6)\a@rgx\a@rgy
            \s@etpospix \a@rgx \a@rgy
            \writeps {\the\t@pixa\space \the\t@pixb\space 
                      \the\t@pixc\space \the\t@pixd\space 
                      \the\x@pix\space \the\y@pix\space cv}}

\def\drawbb {\bsegment
               \drawdim bp
               \setunitscale 0.24
               \linewd 1           
               \writeps {\the\xminpix\space \the\yminpix\space mv}%
               \writeps {\the\xminpix\space \the\ymaxpix\space lv}%
               \writeps {\the\xmaxpix\space \the\ymaxpix\space lv}%
               \writeps {\the\xmaxpix\space \the\yminpix\space lv}%
               \writeps {\the\xminpix\space \the\yminpix\space lv}%
             \esegment}


\def\getpos (#1 #2)#3#4{\g@etargxy #1 #2 {} \\#3#4%
                        \c@heckast #3%
                        \ifa@st
                          \g@etsympix #3\t@pixa
                          \advance \t@pixa by -\x@segoffpix
                          \pixtocoord \t@pixa #3%
                        \fi
                        \c@heckast #4%
                        \ifa@st
                          \g@etsympix #4\t@pixa
                          \advance \t@pixa by -\y@segoffpix
                          \pixtocoord \t@pixa #4%
                        \fi}

\def\getsympos (#1 #2)#3#4{\g@etargxy #1 #2 {} \\#3#4%
                           \c@heckast #3%
                           \ifa@st \else
                             \errmessage {TeXdraw: invalid symbolic coordinate}
                           \fi
                           \c@heckast #4%
                           \ifa@st \else
                             \errmessage {TeXdraw: invalid symbolic coordinate}
                           \fi}

\def\listtopix (#1)#2{\def #2{}
                      \edef\l@ist {#1 }
                      \t@counta=0
                      \loop
                        \expandafter\g@etitem \l@ist \\\a@rgx\l@ist
                        \a@pppix \a@rgx #2%
                        \ifx \l@ist\empty
                          \t@counta=1
                        \fi
                      \ifnum \t@counta=0
                      \repeat}


\def\realmult #1#2#3{\dimen0=#1pt
                     \dimen2=#2\dimen0
                     \edef #3{\expandafter\c@lean\the\dimen2}}

\def\intdiv #1#2#3{\t@counta=#1
                   \t@countb=#2
	           \ifnum \t@countb<0
                      \t@counta=-\t@counta
                      \t@countb=-\t@countb
                   \fi
                   \t@countd=1                    
                   \ifnum \t@counta<0
                      \t@counta=-\t@counta
                      \t@countd=-1
                   \fi
	           \t@countc=\t@counta  \divide \t@countc by \t@countb
                   \t@counte=\t@countc  \multiply \t@counte by \t@countb
                   \advance \t@counta by -\t@counte
	           \t@counte=-1
                   \loop
                     \advance \t@counte by 1
	             \ifnum \t@counte<16
                       \multiply \t@countc by 2           
                       \multiply \t@counta by 2           
                       \ifnum \t@counta<\t@countb \else   
                         \advance \t@countc by 1          
                         \advance \t@counta by -\t@countb 
                       \fi
                   \repeat
	           \divide \t@countb by 2         
	           \ifnum \t@counta<\t@countb     
                     \advance \t@countc by 1
                   \fi
                   \ifnum \t@countd<0             
                     \t@countc=-\t@countc
                   \fi
                   \dimen0=\t@countc sp           
                   \edef #3{\expandafter\c@lean\the\dimen0}}

\outer\def\gnewif #1{\count@\escapechar \escapechar\m@ne
  \expandafter\expandafter\expandafter
   \edef\@if #1{true}{\global\let\noexpand#1=\noexpand\iftrue}%
  \expandafter\expandafter\expandafter
   \edef\@if #1{false}{\global\let\noexpand#1=\noexpand\iffalse}%
  \@if#1{false}\escapechar\count@} 
\def\@if #1#2{\csname\expandafter\if@\string#1#2\endcsname}
{\uccode`1=`i \uccode`2=`f \uppercase{\gdef\if@12{}}} 


\def\coordtopix #1#2{\dimen0=#1\d@dim
                     \dimen2=\d@sc\dimen0
                     \t@counta=\dimen2              
                     \t@countb=\s@ppix
                     \divide \t@countb by 2
                     \ifnum \t@counta<0             
                       \advance \t@counta by -\t@countb
                     \else
                       \advance \t@counta by \t@countb
                     \fi
                     \divide \t@counta by \s@ppix
                     #2=\t@counta}

\def\pixtocoord #1#2{\t@counta=#1%
                     \multiply \t@counta by \s@ppix
                     \dimen0=\d@sc\d@dim
                     \t@countb=\dimen0
                     \intdiv \t@counta \t@countb #2}

\def\pixtodim #1#2{\t@countb=#1%
                   \multiply \t@countb by \s@ppix
                   #2=\t@countb sp\relax}

\def\pixtobp #1#2{\dimen0=\p@sfactor pt
                  \t@counta=\dimen0
                  \multiply \t@counta by #1%
                  \ifnum \t@counta < 0             
                    \advance \t@counta by -32768
                  \else
                    \advance \t@counta by 32768
                  \fi
                  \divide \t@counta by 65536
                  #2=\t@counta}
                  
\newcount\t@counta    \newcount\t@countb   
\newcount\t@countc    \newcount\t@countd
\newcount\t@counte
\newcount\t@pixa      \newcount\t@pixb     
\newcount\t@pixc      \newcount\t@pixd
\let\l@lxbp=\t@pixa   \let\l@lybp=\t@pixb  
\let\u@rxbp=\t@pixc   \let\u@rybp=\t@pixd

\newdimen\t@xpos      \newdimen\t@ypos
\let\l@lxpos=\t@xpos  \let\l@lypos=\t@ypos

\newcount\xminpix      \newcount\xmaxpix
\newcount\yminpix      \newcount\ymaxpix

\newcount\a@lenpix     \newcount\a@widpix

\newcount\x@pix        \newcount\y@pix
\newcount\x@segoffpix  \newcount\y@segoffpix
\newcount\x@savepix    \newcount\y@savepix

\newcount\s@ppix       

\newcount\d@bs

\newcount\t@xdnum
\global\t@xdnum=0

\newdimen\hdrawsize    \newdimen\vdrawsize

\newbox\t@xdbox

\newwrite\drawfile

\newif\ifm@pending
\newif\ifp@ath
\newif\ifa@st
\gnewif \ift@extonly
\gnewif\ifp@osinit

\def\l@paren{(}
\def\a@st{*}

\catcode`\%=12
  \def\p@b {
\catcode`\%=14
\catcode`\{=12  \catcode`\}=12  \catcode`\u=1 \catcode`\v=2
  \def\l@br u{v  \def\r@br u}v
\catcode `\{=1  \catcode`\}=2   \catcode`\u=11 \catcode`\v=11

{\catcode`\p=12 \catcode`\t=12
 \gdef\c@lean #1pt{#1}}

\def\sppix#1/#2 {\dimen0=1#2 \s@ppix=\dimen0
                 \t@counta=#1%
                 \divide \t@counta by 2
                 \advance \s@ppix by \t@counta
                 \divide \s@ppix by #1
                 \t@counta=\s@ppix
                 \multiply \t@counta by 65536       
                 \advance \t@counta by 32891        
                 \divide \t@counta by 65782         
                 \dimen0=\t@counta sp
                 \edef\p@sfactor {\expandafter\c@lean\the\dimen0}}

\def\g@etargxy #1 #2 #3 #4\\#5#6{\def #5{#1}%
                                 \ifx #5\empty
                                   \g@etargxy #2 #3 #4 \\#5#6
                                 \else
                                   \def #6{#2}%
                                   \def\next {#3}%
                                   \ifx \next\empty \else
                                     \errmessage {TeXdraw: invalid coordinate}%
                                   \fi
                                 \fi}

\def\c@heckast #1{\expandafter
                  \c@heckastll #1\\}
\def\c@heckastll #1#2\\{\def\testit {#1}%
                        \ifx \testit\a@st
                          \a@sttrue
                        \else
                          \a@stfalse
                        \fi}

\def\g@etsympix #1#2{\expandafter
                     \ifx \csname #1\endcsname \relax
                       \errmessage {TeXdraw: undefined symbolic coordinate}%
                     \fi
                     #2=\csname #1\endcsname}

\def\s@etcsn #1#2{\expandafter
                  \xdef\csname#1\endcsname {#2}}

\def\g@etitem #1 #2\\#3#4{\edef #4{#2}\edef #3{#1}}
\def\a@pppix #1#2{\edef\next {#1}%
                  \ifx \next\empty \else
                    \coordtopix {#1}\t@pixa
                    \ifx #2\empty
                      \edef #2{\the\t@pixa}%
                    \else
                      \edef #2{#2 \the\t@pixa}%
                    \fi
                  \fi}

\def\s@etpospix #1#2{\coordtopix {#1}\x@pix
                     \advance \x@pix by \x@segoffpix
                     \coordtopix {#2}\y@pix
                     \advance \y@pix by \y@segoffpix
                     \u@pdateminmax \x@pix \y@pix}

\def\r@elpospix #1#2{\coordtopix {#1}\t@pixa
                     \advance \x@pix by \t@pixa
                     \coordtopix {#2}\t@pixa
                     \advance \y@pix by \t@pixa
                     \u@pdateminmax \x@pix \y@pix}

\def\r@elupd #1#2{\t@counta=\x@pix
                  \advance\t@counta by #1%
                  \t@countb=\y@pix
                  \advance\t@countb by #2%
                  \u@pdateminmax \t@counta \t@countb}

\def\u@pdateminmax #1#2{\ifnum #1>\xmaxpix
                          \global\xmaxpix=#1%
                        \fi
                        \ifnum #1<\xminpix
                          \global\xminpix=#1%
                        \fi
                        \ifnum #2>\ymaxpix
                          \global\ymaxpix=#2%
                        \fi
                        \ifnum #2<\yminpix
                          \global\yminpix=#2%
                        \fi}

\def\maxhvpos {\t@pixa=\xmaxpix
               \advance \t@pixa by -\xminpix
               \pixtodim  \t@pixa {\dimen2}%
               \global\hdrawsize=\dimen2
               \t@pixa=\ymaxpix
               \advance \t@pixa by -\yminpix
               \pixtodim \t@pixa {\dimen2}%
               \global\vdrawsize=\dimen2\relax}

\def\t@xdinclude {\pixtobp {-\xminpix}\l@lxbp  \pixtobp {-\yminpix}\l@lybp
                  \ift@extonly \else
                    \includegraphics{\p@sfile\space}%
                  \fi}

\def\s@avemove #1#2{\x@savepix=#1\y@savepix=#2%
                    \m@pendingtrue
                    \ifp@osinit \else
                      \p@osinittrue
                      \global\xminpix=\x@savepix \global\yminpix=\y@savepix
                      \global\xmaxpix=\x@savepix \global\ymaxpix=\y@savepix
                    \fi}

\def\f@lushmove {\p@osinittrue
                 \ifm@pending
                   \writetx {\the\x@savepix\space \the\y@savepix\space mv}%
                   \m@pendingfalse
                   \p@athfalse
                 \fi}

\def\f@lushbs {\loop
                 \ifnum \d@bs>0
                   \writetx {bs}%
                   \global\advance \d@bs by -1
               \repeat}
               
\def\h@move #1#2 #3)#4{\move (#2 #3)%
                       \h@text {#4}}
\def\h@text #1{\pixtodim \x@pix \t@xpos
               \pixtodim \y@pix \t@ypos
               \vbox to 0pt{\normalbaselines
                            \t@stuff
                            \kern -\t@ypos
                            \hbox to 0pt{\l@stuff
                                         \kern \t@xpos
                                         \hbox {#1}%
                                         \kern -\t@xpos
                                         \r@stuff}%
                            \kern \t@ypos
                            \b@stuff\relax}}

\def\r@move td:#1 #2#3 #4)#5{\move (#3 #4)%
                             \r@text td:#1 {#5}}
\def\r@text td:#1 #2{\pixtodim \x@pix \t@xpos
                     \pixtodim \y@pix \t@ypos
                     \vbox to 0pt{\kern -\t@ypos
                                  \hbox to 0pt{\kern \t@xpos
                                               \rottxt{#1}{#2}%
                                               \hss}%
                                  \vss}}

\def\rottxt #1#2{
                 \z@sb{#2}%
}
\def\z@sb #1{\vbox to 0pt{\normalbaselines
                          \t@stuff
                          \hbox to 0pt{\l@stuff
                                       \hbox {#1}%
                                       \r@stuff}%
                          \b@stuff}}

\def\t@exdrawdef {\sppix 300/in            
                  \drawdim in              
                  \edef\u@nitsc {1}
                  \setsegscale 1           
                  \arrowheadsize l:0.16 w:0.08
                  \arrowheadtype t:T
                  \textref h:L v:B }


\def\writeps #1{\f@lushbs
                \f@lushmove
                \p@athtrue
                \writetx {#1}}
\def\writetx #1{\ift@extonly
                  \t@extonlyfalse
                  \t@dropen
                \fi
                \w@rps {#1}}
\def\w@rps #1{\immediate\write\drawfile {#1}}

\def\t@dropen {%
  \global\advance \t@xdnum by 1
  \ifnum \t@xdnum<10
    \xdef\p@sfile {\jobname.ps\the\t@xdnum}
  \else
    \xdef\p@sfile {\jobname.p\the\t@xdnum}
  \fi
  \immediate\openout\drawfile=\p@sfile
  \w@rps {\p@b PS-Adobe-3.0 EPSF-3.0}%
  \w@rps {\p@p BoundingBox: (atend)}%
  \w@rps {\p@p Title: TeXdraw drawing: \p@sfile}%
  \w@rps {\p@p Pages: 1}%
  \w@rps {\p@p Creator: \TeXdrawId}%
  \w@rps {\p@p CreationDate: \the\year/\the\month/\the\day}%
  \w@rps {50 dict begin}%
  \w@rps {/mv {stroke moveto} def}%
  \w@rps {/lv {lineto} def}%
  \w@rps {/st {currentpoint stroke moveto} def}%
  \w@rps {/sl {st setlinewidth} def}%
  \w@rps {/sd {st 0 setdash} def}%
  \w@rps {/sg {st setgray} def}%
  \w@rps {/bs {gsave} def /es {stroke grestore} def}%
  \w@rps {/cv {curveto} def}%
  \w@rps {/cr \l@br gsave /rad exch def currentpoint newpath rad 0 360 arc}%
  \w@rps { stroke grestore\r@br\space def}%
  \w@rps {/fc \l@br gsave /rad exch def setgray currentpoint newpath}%
  \w@rps { rad 0 360 arc fill grestore\r@br\space def}%
  \w@rps {/ar {gsave currentpoint newpath 5 2 roll arc stroke grestore} def}%
  \w@rps {/el \l@br gsave /rady exch def /radx exch def}%
  \w@rps { /svm matrix currentmatrix def currentpoint translate}%
  \w@rps { radx rady scale newpath 0 0 1 0 360 arc}%
  \w@rps { svm setmatrix stroke grestore\r@br\space def}%
  \w@rps {/fl \l@br gsave closepath setgray fill grestore}%
  \w@rps { currentpoint newpath moveto\r@br\space def}%
  \w@rps {/fp \l@br gsave closepath setgray fill grestore}%
  \w@rps { currentpoint stroke moveto\r@br\space def}%
  \w@rps {/av \l@br /hhwid exch 2 div def /hlen exch def}%
  \w@rps { /ah exch def /tipy exch def /tipx exch def}%
  \w@rps { currentpoint /taily exch def /tailx exch def}%
  \w@rps { /dx tipx tailx sub def /dy tipy taily sub def}%
  \w@rps { /alen dx dx mul dy dy mul add sqrt def}%
  \w@rps { /blen alen hlen sub def}%
  \w@rps { gsave tailx taily translate dy dx atan rotate}%
  \w@rps { (V) ah ne {blen 0 gt {blen 0 lineto} if} {alen 0 lineto} ifelse}%
  \w@rps { stroke blen hhwid neg moveto alen 0 lineto blen hhwid lineto}%
  \w@rps { (T) ah eq {closepath} if}%
  \w@rps { (W) ah eq {gsave 1 setgray fill grestore closepath} if}%
  \w@rps { (F) ah eq {fill} {stroke} ifelse}%
  \w@rps { grestore tipx tipy moveto\r@br\space def}%
  \w@rps {\p@sfactor\space \p@sfactor\space scale}%
  \w@rps {1 setlinecap 1 setlinejoin}%
  \w@rps {3 setlinewidth [] 0 setdash}%
  \w@rps {0 0 moveto}%
}

\def\t@drclose {%
  \w@rps {stroke end showpage}%
  \w@rps {\p@p Trailer:}%
  \pixtobp \xminpix \l@lxbp  \pixtobp \yminpix \l@lybp
  \pixtobp \xmaxpix \u@rxbp  \pixtobp \ymaxpix \u@rybp
  \w@rps {\p@p BoundingBox: \the\l@lxbp\space \the\l@lybp\space
                            \the\u@rxbp\space \the\u@rybp}%
  \w@rps {\p@p EOF}%
  \closeout\drawfile
}

\catcode`\@=\catamp